\documentclass[seceqn,secthm]{elsart}
\usepackage{amssymb}\usepackage{amsmath}
\usepackage[all, knot]{xy}
\xyoption{arc} \xyoption{color}\xyoption{line}
\usepackage{mathrsfs}

\renewcommand{\leq}{\leqslant}
\renewcommand{\geq}{\geqslant}

\newcommand{\supp}{\operatorname{supp}}
\newcommand{\diag}{\operatorname{diag}}

\journal{arXiv:}

\begin{document}

\begin{frontmatter}
\title{Exponential Decay of Eigenfunctions and Accumulation of Eigenvalues  on Manifolds with Axial Analytic Asymptotically Cylindrical Ends\thanksref{AKA}}
\author{Victor Kalvin}
\thanks[AKA]{This work was funded by grant N108898 awarded by the Academy of Finland.}
\ead{vkalvin\,@\,gmail.com}
\date{}

\begin{abstract} In this paper we continue our study of the Laplacian on manifolds with axial analytic asymptotically cylindrical ends initiated in~arXiv:1003.2538. By using the complex scaling method and the Phragm\'{e}n-Lindel\"{o}f principle we prove exponential decay of the eigenfunctions corresponding to the non-threshold eigenvalues of the Laplacian on functions. In the case of a manifold with (non-compact) boundary it is either the Dirichlet Laplacian or the Neumann Laplacian. We show that the rate of exponential decay of an eigenfunction is prescribed by the distance from the corresponding eigenvalue to the next threshold. Under our assumptions on the behaviour of the metric at infinity accumulation of isolated and embedded eigenvalues occur. The results on decay of eigenfunctions combined with the compactness argument due to Perry imply that the eigenvalues can accumulate only at thresholds and only from below. The eigenvalues are of finite multiplicity.
\end{abstract}

\end{frontmatter}
\section{Introduction}
Consider a manifold  $(\mathcal M,\mathsf g)$ with an asymptotically cylindrical end. This means that $\mathcal M$ is a smooth non-compact manifold of the form $\mathcal M_c\cup (\Bbb R_+\times\Omega)$, where $\mathcal M_c$ is a compact manifold, and  $\Bbb R_+\times \Omega$ is the Cartesian product of the positive semi-axis $\Bbb R_+$ and a compact manifold $\Omega$, see Fig.~\ref{manif0} and Fig.~\ref{manif}. Furthermore, the  metric $\mathsf g$ asymptotically approaches at infinity  the product metric $dx\otimes dx+\mathfrak h$ on the semi-cylinder $\Bbb R_+\times\Omega$, where $\mathfrak h$ is a metric on $\Omega$. Traditionally, one studies the Laplacian on  manifolds with asymptotically cylindrical ends under
more~\cite{chr ST,Guillope,Melrose,MelroseScat} or less~\cite{Edward,FroHislop,Mueller,IKL} restrictive assumptions on the rate of convergence of the metric $\mathsf g$ to the product metric $dx\otimes dx+\mathfrak h$ at infinity.
We do not make this kind of assumptions. Instead,  in~\cite{KalvinII} and in this paper we consider  manifolds with axial analytic asymptotically cylindrical ends. On these manifolds the metric $\mathsf g$ extends by analyticity to a conical neighborhood of the axis $\Bbb R_+$ of the semi-cylinder $\Bbb R_+\times \Omega$, and the continuation tends at infinity to the analytic continuation of $dx\otimes dx+\mathfrak h$; for precise definitions see Section~\ref{s2}. Due to these properties of $\mathsf g$ the end $(\Bbb R_+\times \Omega,\mathsf g\!\upharpoonright_{\Bbb R_+\times \Omega})$ is said to be axial analytic.
On manifolds with axial analytic asymptotically cylindrical ends the metric $\mathsf g$ converges at infinity to  $dx\otimes dx+\mathfrak h$ with arbitrarily slow rate.

As is known~\cite{CZ,P}, the eigenvalues  of the Laplacian on manifolds with cylindrical ends have no finite points of accumulation.  In the cylindrical ends we have $\mathsf g\!\upharpoonright_{\Bbb R_+\times\Omega}=dx\otimes dx +\mathfrak h$, and exponential decay of the non-threshold eigenfunctions can be easily seen by separation of variables. Let us also note that finite points of accumulation do not occur on manifolds with asymptotically cylindrical ends, provided that we admit only  exponential convergence of the metric $\mathsf g$ to $dx\otimes dx +\mathfrak h$ at infinity~\cite{Melrose}. In this case exponential decay of the non-threshold eigenfunctions is a consequence of the asymptotic theory, see e.g.~\cite{KozlovMaz`ya,KozlovMazyaRossmann,Melrose} and references therein. Once we allow for sufficiently slow convergence of $\mathsf g$ to $dx\otimes dx +\mathfrak h$ at infinity, the situation changes and eigenvalues may  accumulate at finite distances. Under an assumption on the rate of convergence of the metric at infinity it is possible to prove the Mourre estimates~\cite{FroHislop}. In particular, these estimates imply that the eigenvalues can accumulate only at thresholds.  Power decay of eigenfunctions and eigenvalue accumulation for the Laplacian on functions, for the Dirichlet Laplacian, and for the Neumann Laplacian were studied in~\cite{Edward}, where it is assumed that the metric allows for separation of variables in the ends and satisfies an assumption on the rate of convergence at infinity.

 In~\cite{KalvinII}
we developed an approach to the complex scaling on manifolds  with  axial analytic asymptotically cylindrical ends and
established a variant of the Aguilar-Balslev-Combes theorem for the Laplacian $\Delta$ on functions.
In particular, we proved that the Laplacian has no singular continuous spectrum, all non-threshold eigenvalues are of finite multiplicity, and the eigenvalues of $\Delta$ can accumulate only at thresholds.  In this paper we continue to study the Laplacian and prove the following: 1)~Any non-threshold eigenfunction of $\Delta$ decays at infinity with some exponential rate prescribed by the distance from the corresponding eigenvalue to the next threshold of $\Delta$; 2)~The eigenvalues of $\Delta$ are of finite multiplicity and can accumulate at thresholds only from below. Besides, we aim to show a certain similarity between methods and results of the theory of $N$-body Schr\"{o}dinger operators in $\Bbb R^n$
and the analysis on manifolds with axial analytic asymptotically cylindrical ends. It is interesting to note that there is also a connection between the theory of Schr\"{o}dinger operators and the analysis on symmetric spaces~\cite{MV1,MV2}. As in~\cite{KalvinII} we consider three generic cases:
1)~$\Delta$ is the Laplacian on a manifold without boundary;
2)~$\Delta$ is the Dirichlet Laplacian on a manifold with noncompact boundary;
3)~$\Delta$ is the Neumann Laplacian on a manifold with noncompact boundary.
We give examples of manifolds with  axial analytic asymptotically cylindrical ends demonstrating that the eigenvalues of $\Delta$ may indeed accumulate at thresholds.

As is typically the case, it is much easier to prove exponential decay of the eigenfunctions corresponding to the isolated eigenvalues.  In fact, this can be done by methods of the asymptotic theory~\cite{KozlovMaz`ya,KozlovMazyaRossmann,MP2}, which work nicely on manifolds with axial analytic asymptotically cylindrical ends  (see also~\cite{Kalvine} and~\cite[Appendix]{PhD}).  Or, equivalently, one can use the dilation analytic techniques similar to those in the theory of Schr\"{o}dinger operators e.g.~\cite[Chapter XII.11]{Simon Reed iv}. However,  $\Delta$ is a non-negative operator,  and therefore  only the Dirichlet Laplacian may have isolated eigenvalues below its absolutely continuous spectrum $\sigma_{ac}(\Delta)=[\nu,\infty)$, $\nu >0$. All eigenvalues of the Neumann Laplacian and of  the Laplacian on  a manifold $(\mathcal M,\mathsf g)$ without boundary  are embedded into the  spectrum $\sigma_{ac}(\Delta)=[0,\infty)$.

As is well-known, for Schr\"{o}dinger operators with dilation analytic potentials it is also possible to prove exponential decay of the eigenfunctions corresponding to the  non-threshold embedded eigenvalues, see e.g.~\cite[Chapter XII.11]{Simon Reed iv}. It turns out that similar methods, based on the complex scaling and the Phragm\'{e}n-Lindel\"{o}f principle, can be applied on manifolds with axial analytic asymptotically  cylindrical ends.
This allows us to prove that every eigenfunction  corresponding to a non-threshold eigenvalue of the Laplacian is of some exponential decay at infinity. This fact plays a crucial role throughout the paper. Note that due to arbitrarily slow convergence of $\mathsf g$ to $dx\otimes dx +\mathfrak h$ at infinity the asymptotic theory does not give any information on decay of eigenfunctions corresponding to an embedded eigenvalue, see e.g.~\cite{Kalvine,KozlovMaz`ya,KozlovMazyaRossmann,MP2,Plam}. Nonetheless, one can employ the asymptotic theory in order to refine the rate of exponential decay of eigenfunctions. As a result we conclude that the rate of exponential decay of a non-threshold eigenfunction is prescribed by the distance from the corresponding eigenvalue to the next threshold of the Laplacian. Similar results for  $N$-body Schr\"{o}dinger operators can be found in~\cite{FH}, see also references therein.

Finally, we study accumulation of eigenvalues. Here we  combine our results on decay of eigenfunctions with the compactness argument due to Perry~\cite{Perry}. Let us remark that an attempt to describe  accumulation of eigenvalues  for general elliptic selfadjoint problems in domains with cylindrical ends was made in~\cite{KNP}. However, as it was observed later~\cite{Plam}, the results on accumulation of eigenvalues announced in~\cite{KNP} are valid only under an additional assumption on exponential decay of eigenfunctions;  see~\cite{Kalvine} for the proof of other results announced in~\cite{KNP}. Despite that after simple modifications the approach~\cite{Kalvine,KNP,Plam} is capable to prove some results of this paper,  we prefer to rely on methods similar to those we meet in the theory of Schr\"{o}dinger operators~\cite{Cycon,FH,Hunziker,Perry, Simon Reed iv}. First, because this demonstrates a certain similarity between the theory of Schr\"{o}dinger operators
and the analysis on manifolds with axial analytic asymptotically cylindrical ends. Secondly, because these methods are simpler: in contrast to~\cite{Kalvine,KNP,Plam} they do not require from the reader an extensive background in the analytic Fredholm theory~\cite{GS} nor any prior knowledge of methods and results of the asymptotic theory~\cite{KozlovMaz`ya,KozlovMazyaRossmann,MP2}. Nonetheless, throughout the paper we make corresponding remarks every time one assertion or another can equivalently be obtained by methods of~\cite{Kalvine,KNP,KozlovMaz`ya,KozlovMazyaRossmann,MP2,Plam}.

The structure of this paper is as follows. In Section~\ref{s2} we introduce manifolds with asymptotically cylindrical and axial analytic asymptotically cylindrical ends. Section~\ref{s3} presents a summary of main results of this paper and two illustrative examples. In Section~\ref{sCS} we give a broad overview of our approach to the complex scaling~\cite{KalvinII}. In Section~\ref{s5} we study the quadratic form and localize the essential spectrum of the Laplacian deformed by means of the complex scaling and conjugated with an exponent.
Then in Section~\ref{s6} we prove that all non-threshold eigenfunctions of the Laplacian are of some exponential decay at infinity. Finally, in Section~\ref{s7} we show that the rate of eigenfunction decay is prescribed by the distance from the corresponding eigenvalue to the next threshold, and study accumulation of eigenvalues.

\section{Manifolds with axial analytic asymptotically cylindrical ends}\label{s2}
\begin{figure}[h]
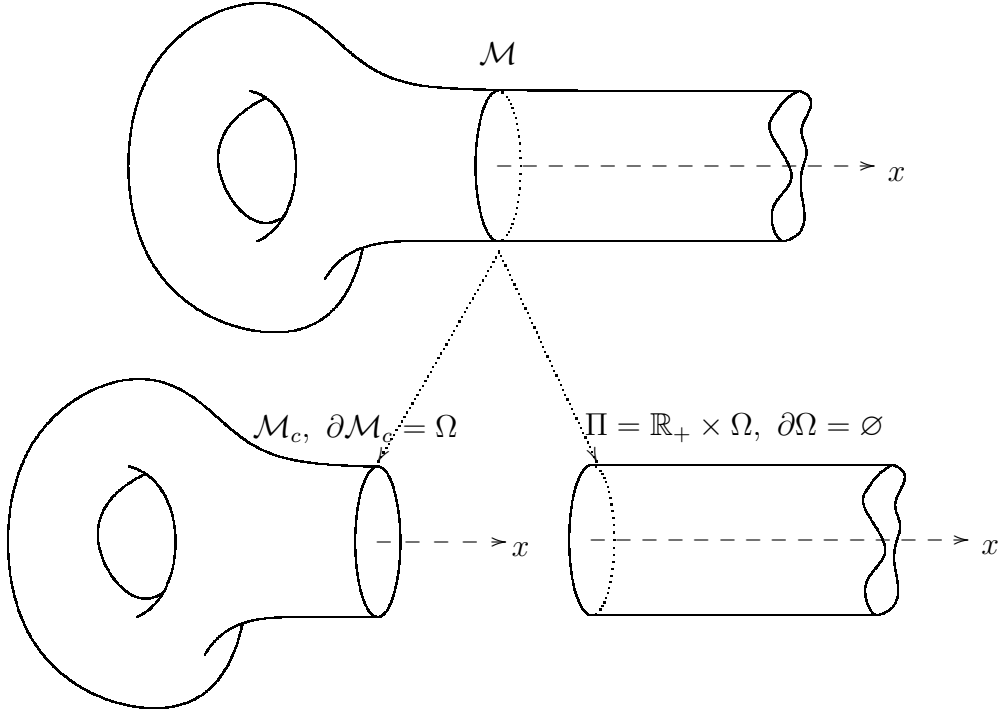

\[\xy (33,50)*{\xy
(40,10);(-18,-11)**\crv{(0,10)&(-13,10)&(-20,15)&(-30,25)&(-50,15)&(-50,-15)&(-30,-25)&(-20,-20)};
(38,-10);(-23,-15)**\crv{(0,-10)&(-10,-10)&(-20,-10)};
(-32,-10); (-33,10)**\crv{(-25,-7)&(-25,8)};
(-29,-7); (-31,9)**\crv{(-32,-9)&(-39,-2)&(-37,6)};
(40,10); (38,-10)**\crv{(43,9)&(39,5)&(42,2)&(39,-2)&(42,-9)};
(40,10); (38,-10)**\crv{(38,9)&(37,7)&(35,5)&(41,-2)&(35,-5)&(37,-10)};
{\ar@{-->} (0,0)*{}; (50,0)*{}}; (53,-1)*{x};
(0,0)*\ellipse(3,10){.};
(0,0)*\ellipse(3,10)__,=:a(-180){-};
(0,15)*{\mathcal M};
\endxy};
(0,0)*{\xy (0,10); (-18,-11)**\crv{(0,10)&(-13,10)&(-20,15)&(-30,25)&(-50,15)&(-50,-15)&(-30,-25)&(-20,-20)};
(0,-10); (-23,-15)**\crv{(0,-10)&(-10,-10)&(-20,-10)};
(-32,-10); (-33,10)**\crv{(-25,-7)&(-25,8)};
(-29,-7); (-31,9)**\crv{(-32,-9)&(-39,-2)&(-37,6)};
(-3,15)*{\mathcal M_c,\ \partial\mathcal M_c=\Omega};
{\ar@{-->} (0,0)*{}; (17,0)*{}}; (19,-1)*{x};
(0,0)*\ellipse(3,10){-};
(0,0)*\ellipse(3,10)__,=:a(-180){-};
\endxy};
(70,4)*{\xy (0,10); (40,10)**\dir{-};(0,-10); (38,-10)**\dir{-};
(0,0)*\ellipse(3,10){.};
(0,0)*\ellipse(3,10)__,=:a(-180){-};
(19,15)*{\Pi=\Bbb R_+\times\Omega,\ \partial\Omega=\varnothing};(0,15)*{\vphantom{\mathcal M_c}};
(40,10); (38,-10)**\crv{(43,9)&(39,5)&(42,2)&(39,-2)&(42,-9)};
(40,10); (38,-10)**\crv{(38,9)&(37,7)&(35,5)&(41,-2)&(35,-5)&(37,-10)};
{\ar@{-->} (0,0)*{}; (50,0)*{}}; (53,-1)*{x};
(0,0)*\ellipse(3,10){.};
(0,0)*\ellipse(3,10)__,=:a(-180){-};
\endxy};
{\ar@{.>} (31,39)*{};(15,11)*{}};{\ar@{.>} (31,39)*{};(44,11)*{}};
\endxy\]
\caption{Representation $\mathcal M=\mathcal M_{c}\cup\Pi$ of a manifold $\mathcal M$ without boundary.}\label{manif0}
\end{figure}

Let $\Omega$ be a smooth compact  $n$-dimensional manifold with smooth boundary $\partial\Omega$ or without it.
Denote by $\Pi$ the semi-cylinder $\Bbb R_+\times \Omega$, where $\Bbb R_+$ is the  positive  semi-axis, and $\times$ stands for the Cartesian product. Consider  a  smooth oriented connected $n+1$-dimensional  manifold $\mathcal M$   representable in the form $\mathcal M=\mathcal M_c\cup
\Pi$, where $\mathcal M_c$ is a smooth compact manifold with  boundary, cf. Fig.~\ref{manif0} and Fig.~\ref{manif}. We do not  consider the case of a manifold $\mathcal M$ with compact boundary $\partial \mathcal M$ as it can be treated similarly to the case of a manifold without boundary.
 Namely, we assume that $\partial\mathcal M=\varnothing$  in the case $\partial\Omega=\varnothing$.
\begin{figure}[h]
\[\xy (33,50)*{\xy
(40,10);(-18,-11)**\crv{(0,10)&(-13,10)&(-20,15)&(-30,25)&(-50,15)&(-50,-15)&(-30,-25)&(-20,-20)};
(38,-10);(-23,-15)**\crv{(0,-10)&(-10,-10)&(-20,-10)};
(-32,-10); (-33,10)**\crv{(-25,-7)&(-25,8)};
(-29,-7); (-31,9)**\crv{(-32,-9)&(-39,-2)&(-37,6)};
(40,10); (38,-10)**\crv{(43,9)&(39,5)&(42,2)&(39,-2)&(42,-9)};
(40,10); (38,-10)**\crv{(38,9)&(37,7)&(35,5)&(41,-2)&(35,-5)&(37,-10)};
{\ar@{-->} (0,0)*{}; (50,0)*{}}; (53,-1)*{x};
(0,0)*\ellipse(3,10){.};
(0,0)*\ellipse(3,10)__,=:a(-180){-};
(0,15)*{\mathcal M};
\endxy};
(0,0)*{\xy (8,10); (-18,-11)**\crv{(0,10)&(-13,10)&(-20,15)&(-30,25)&(-50,15)&(-50,-15)&(-30,-25)&(-20,-20)};
(8,-10); (-23,-15)**\crv{(0,-10)&(-10,-10)&(-20,-10)};
(-32,-10); (-33,10)**\crv{(-25,-7)&(-25,8)};
(-29,-7); (-31,9)**\crv{(-32,-9)&(-39,-2)&(-37,6)};
(-5,15)*{\mathcal M_c};
(10,-7); (10,6)**\crv{(10,-6)&(14,-2)&(15,1)&(10,4)&(10,7)};
(8,10); (13,2)**\crv{(12,10)};
(8,-10); (13,-3)**\crv{(12,-10)};
{\ar@{-->} (0,0)*{}; (17,0)*{}}; (19,-1)*{x};
(0,0)*\ellipse(3,10){.};
(0,0)*\ellipse(3,10)__,=:a(-180){-};
\endxy};
(70,4)*{\xy (0,10); (40,10)**\dir{-};(0,-10); (38,-10)**\dir{-};
(0,0)*\ellipse(3,10){.};
(0,0)*\ellipse(3,10)__,=:a(-180){-};
(21,15)*{\Pi=\Bbb R_+\times\Omega,\ \partial\Omega\neq\varnothing};(0,15)*{\vphantom{\mathcal M_c}};
(40,10); (38,-10)**\crv{(43,9)&(39,5)&(42,2)&(39,-2)&(42,-9)};
(40,10); (38,-10)**\crv{(38,9)&(37,7)&(35,5)&(41,-2)&(35,-5)&(37,-10)};
{\ar@{-->} (0,0)*{}; (50,0)*{}}; (53,-1)*{x};
(0,0)*\ellipse(3,10){.};
(0,0)*\ellipse(3,10)__,=:a(-180){-};
\endxy};
{\ar@{.>} (31,39)*{};(15,11)*{}};{\ar@{.>} (31,39)*{};(44,11)*{}};
\endxy\]
\caption{Representation $\mathcal M=\mathcal M_{c}\cup\Pi$ of a  manifold $\mathcal M$ with boundary.}\label{manif}
\end{figure}

Let $\mathsf g\in C^\infty \mathrm T^*\mathcal M^{\otimes 2}$ be a  Riemannian metric on $\mathcal M$.
We identify the cotangent bundle $\mathrm T^*\Pi$ with the tensor product $\mathrm T^*\Bbb R_+\otimes\mathrm T^*\Omega$   via the natural isomorphism induced by the product structure on $\Pi$. This together with the trivialization $\mathrm T^*\Bbb R_+=\{(x,a\, dx): x\in\Bbb R_+,a\in \Bbb R\}$ implies that any metric   $\mathsf g$ can be represented on $\Pi$ in the form
\begin{equation}\label{split}
\mathsf g\!\upharpoonright_{\Pi}=\mathfrak g_0 dx\otimes dx+2 \mathfrak g_1\otimes dx +\mathfrak g_2,
\quad\mathfrak g_k(x)\in C^\infty\mathrm T^*\Omega^{\otimes k}.
\end{equation}

Denote by  $\Bbb C\mathrm T^*\Omega^{\otimes k}$  the tensor power of the complexified cotangent bundle $\Bbb C\mathrm T^*\Omega$ with the fibers $\Bbb C\mathrm T_{\mathrm y}^*\Omega=\mathrm T_{\mathrm y}^*\Omega\otimes\Bbb C$. In what follows $C^m$ stands for sections of complexified bundles, e.g. we write $C^\infty \mathrm T^*\Omega^{\otimes k}$ and $C^1 \mathrm T^*\Omega^{\otimes k}$ instead of $C^\infty\Bbb C \mathrm T^*\Omega^{\otimes k}$ and $C^1\Bbb C\mathrm T^*\Omega^{\otimes k}$.
We equip the space $C^1\mathrm T^*\Omega^{\otimes k}$ with the norm
\begin{equation}\label{Enorm}
\|\cdot\|_{\mathfrak e}=\max_{\mathrm y\in\Omega}\bigl(|\cdot|_{\mathfrak e}(\mathrm y)+|D\cdot|_{\mathfrak e}(\mathrm y)\bigr),
\end{equation}
where $\mathfrak e$ is a Riemannian metric on $\Omega$, $|\cdot|_{\mathfrak e}(\mathrm y)$ is the norm induced by  $\mathfrak e$ in the fiber $\Bbb C\mathrm T_{\mathrm y}^*\Omega^{\otimes k}$, and
$ D: C^1\mathrm T^*\Omega^{\otimes k}\to C^0\mathrm T^*\Omega^{\otimes k+1}$
is the Levi-Civita connection on the manifold $(\Omega,\mathfrak e)$.

\begin{defn}\label{ACE} We say that $(\mathcal M,\mathsf g)$ is a manifold with an axial analytic asymptotically cylindrical end $(\Pi,\mathsf g\!\upharpoonright_{\Pi})$, if the following conditions hold:
\begin{itemize}
\item[i.]  The functions   $ x\mapsto \mathfrak g_k(x)\in C^\infty\mathrm T^*\Omega^{\otimes k}$ in~\eqref{split}   extend  by analyticity in $x$ from $\Bbb R_+$ to the sector $\Bbb S_\alpha=\{z\in\Bbb C:|\arg z|<\alpha\}$ with some  $\alpha>0$.
\item[ii.] As $z$ tends to infinity in $\Bbb S_\alpha$ the function $\mathfrak g_0(z)$ uniformly converges to $1$ in the norm of $C^1(\Omega)$, the tensor field $\mathfrak g_1(z)$ uniformly converges to zero in the norm of $C^1\mathrm T^*\Omega$, and the tensor field $\mathfrak g_2(z)$ uniformly converges to a Riemannian metric $\mathfrak h$ on $\Omega$ in the norm of $C^1\mathrm T^*\Omega^{\otimes 2}$.
\end{itemize}
\end{defn}

In this paper we are mainly concerned in manifolds  with axial analytic asymptotically cylindrical ends.
However, some of our results are valid without any  assumptions on the axial analytic regularity of the  metric $\mathsf g\!\upharpoonright_\Pi$. In those cases  we consider general manifolds  with asymptotically cylindrical ends in the sense of the following definition.
\begin{defn}\label{ace} We say that $(\mathcal M,\mathsf g)$ is a manifold with an asymptotically cylindrical end $(\Pi,\mathsf g\!\upharpoonright_{\Pi})$, if
for some Riemannian metric $\mathfrak h$  on $\Omega$ we have
$$
\|\mathfrak g_0(x)-1\|_{\mathfrak e}+\|\mathfrak g_1(x)\|_{\mathfrak e}+\|\mathfrak g_2(x)-\mathfrak h\|_{\mathfrak e}\to 0\text{ as } x\to +\infty,
$$
and also
$$\max_{\mathrm y\in\Omega}\sum_{k=0}^2|\partial_x \mathfrak g_k(x)|_{\mathfrak e}(\mathrm y)\to 0\text{ as } x\to+\infty,$$
 where $\mathfrak g_k$ are the coefficients in~\eqref{split}, and $\partial_x=d/dx$ is the real derivative.
\end{defn}
Note that Definitions~\ref{ace} and~\ref{ACE} are independent of the metric $\mathfrak e$ on $\Omega$. It is a consequence of the Cauchy inequalities that any manifold with an axial analytic asymptotically cylindrical end is also a manifold with an asymptotically cylindrical end.

\section{Summary of main results}\label{s3}
Consider a manifold $(\mathcal M,\mathsf g)$  with an asymptotically cylindrical end $(\Pi,\mathsf g\!\upharpoonright_{\Pi})$. We introduce the Hilbert space $L^2(\mathcal M)$ as the completion of the set $C_c^\infty(\mathcal M)$ with respect to the norm $\|\cdot\|=\sqrt{(\cdot,\cdot)}$, where $(\cdot,\cdot)$ is the global inner product on $(\mathcal M,\mathsf g)$.
Let $\Delta$ be the Laplacian on $(\mathcal M, \mathsf g)$ initially defined on a  core $\mathbf C(\Delta)$. In the case
$\partial\mathcal M=\varnothing$  we take  $\mathbf C(\Delta)\equiv C^\infty_c(\mathcal M)$, while in the case $\partial\mathcal M\neq \varnothing$ the core $\mathbf C({\Delta})$ of the Neumann (resp. Dirichlet) Laplacian $\Delta$  consists of the functions $u \in
C^\infty_c({\mathcal M})$ satisfying the Neumann boundary condition  ${{\partial}_\nu} u=0$ (resp. the Dirichlet boundary condition $u\!\upharpoonright_{\partial \mathcal M}=0$).

Let $(\Omega,\mathfrak h)$ be the same compact Riemannian manifold
as in Definition~\ref{ACE}. Recall that we exclude from consideration the case of a manifold $\mathcal M$ with compact boundary $\partial\mathcal M$, i.e. the equalities $\partial\mathcal M=\varnothing$ and $\partial\Omega=\varnothing$ can hold only simultaneously.  If $\partial
\mathcal M\neq \varnothing$ and  $\Delta$ is the Dirichlet (resp. Neumann) Laplacian on $(\mathcal M,\mathsf
g)$, then  by $\Delta_\Omega$ we denote the Dirichlet (resp. Neumann)  Laplacian on  $(\Omega,\mathfrak h)$. If $\partial\mathcal M=\varnothing$, then $\Delta_\Omega$ is the Laplacian on the manifold $(\Omega,\mathfrak h)$ without boundary. Let $L^2(\Omega)$ be the Hilbert space of all square summable functions on $(\Omega,\mathfrak h)$. As is well-known, the spectrum of the operator
$\Delta_\Omega$ in $L^2(\Omega)$ consists of  infinitely many  nonnegative isolated eigenvalues. Let $\nu_1<\nu_2<\dots$ be the distinct
eigenvalues of $\Delta_\Omega$. By definition
$\{\nu_j\}_{j=1}^\infty$ is the set of thresholds of the Laplacian $\Delta$ on $(\mathcal M,\mathsf g)$. Let us
stress that the  thresholds of the Dirichlet and the Neumann Laplacians  are different.

The  main results of this paper are listed in the next theorem.
\begin{thm}\label{main} Let $\Delta$ be the Laplacian on a manifold $(\mathcal M,\mathsf g)$ with an axial analytic asymptotically cylindrical end. Then the following assertions are valid.
\begin{itemize}
\item[1.] The eigenvalues of the selfadjoint operator $\Delta$ in $L^2(\mathcal M)$ are of finite multiplicity and   can accumulate only at the thresholds $\nu_1,\nu_2,\dots$, and only from below.
\item[2.] Any eigenfunction $\Psi$  corresponding to a non-threshold  eigenvalue $\mu$ of $\Delta$ meets the estimate
    \begin{equation}\label{pointwise}
\|\Psi(x)\|_{ L^2(\Omega)}\leq C e^{\gamma x}\text{ as } x\to+\infty
\end{equation}
 with any negative $\gamma>-\min_{j:\nu_j>\mu}{\sqrt{\nu_j-\mu}}$ and an independent of $x$ constant $C$. Here $x\in\Bbb R_+$ is the axial coordinate of the end $\Pi$, see Fig.~\ref{manif0} and Fig.~\ref{manif}, and $\min_{j:\nu_j>\mu}(\nu_j-\mu)$ is the distance from $\mu$ to the next threshold of $\Delta$.
\end{itemize}
\end{thm}
Similar results for $N$-body Schr\"{o}dinger operators can be found e.g. in~\cite{FH,Perry,Simon Reed iv}.
We complete this section with examples of manifolds with axial analytic asymptotically cylindrical ends for which the eigenvalues of the Laplacian accumulate at thresholds.

Consider a smooth compact $n$-dimensional Riemannian  manifold  $(\Omega',\mathfrak h$) with smooth boundary or without it. Let the infinite cylinder $\mathcal M=\Bbb R\times\Omega'$ be endowed with the metric $\mathsf g=d x\otimes dx+f(x)^{4/n} \mathfrak h$, where $f$ is a smooth positive function on $\Bbb R$, such that  $f(x)^{4/n}=1+|x|^{-\delta}$ for   $|x|\geq c>0$ and  $\delta\in(0,2]$. Then $(\mathcal M,\mathsf g)$ can be viewed as a manifold with the axial analytic  asymptotically cylindrical end $(\Bbb R_+\times\Omega, d x\otimes dx+f(x+c)^{4/n} \mathfrak h)$, where $\Omega$ consists of two copies of $\Omega'$.  Let $\{\sigma_k\}_{k=1}^\infty$ be the  eigenvalues of $\Delta_{\Omega'}$ listed with multiplicity. Separation of variables~\cite{CZ} shows that the Laplacian $\Delta=\frac{1}{f} (-\partial_x^2+f''/f+f^{-4/n}\Delta_{\Omega'})f$
on $(\mathcal M,\mathsf g)$ (it is the Dirichlet Laplacian  in the case $\partial\mathcal M\neq\varnothing$) is unitary equivalent to the direct sum of the unbounded operators
$-\partial_x^2+f''/f+f^{-4/n}\sigma_k$
 acting in $L^2(\Bbb R)$. Therefore $\mu$ is an eigenvalue of $\Delta$, if and only if $\mu-\sigma_k$ is an eigenvalue of the Schr\"{o}dinger operator $-\partial_x^2+V_k$ in $L^2(\Bbb R)$ for some $k$, where $V_k=f''/f+(f^{-4/n}-1)\sigma_k$ is the potential.  The minimax principle implies that for all sufficiently large $\sigma_k$ the discrete eigenvalues of  $-\partial_x^2+V_k$  accumulate at zero from below, cf.~\cite[Theorem XIII.6]{Simon Reed iv}. Thus the embedded eigenvalues of $\Delta$  accumulate at every sufficiently large threshold $\nu_j=\sigma_k$.

 Let $f\in C^\infty(\Bbb R)$, $f(s)> 0$, and $f(s)=1+5|s|^{-\delta}$ for   $|s|\geq c>0$ and  $\delta\in(0,2]$.   The domain   $\mathcal G=\{(s,t)\in \Bbb R^2: |t|\leq f(s)\}$ can be viewed as a manifold $(\mathcal M,\mathsf g)$ with an axial analytic asymptotically cylindrical end.  Indeed, we can set $\Pi=\Bbb R_+\times \{[-1,1]\cup[-1,1]\}$ and define $\mathsf g\!\upharpoonright_\Pi$ as the pullback of the Euclidean metric by the diffeomorphism $(\pm s,t)=(x+c,f(x+c) y)$ mapping $\Pi$ onto the asymptotic semi-strips  $\{(s,t)\in \mathcal G:  \pm s>c\}$; here $y\in[-1,1]$ is the local coordinate on $\Omega=[-1,1]\cup[-1,1]$. Due to the axial symmetry of $\mathcal G$ it is possible to prove by the minimax principle that the embedded eigenvalues of the Neumann Laplacian on $(\mathcal M,\mathsf g)$ accumulate at the first non-zero threshold $\nu_2=\pi^2/4$, while the isolated and embedded eigenvalues of the Dirichlet Laplacian  on $(\mathcal M,\mathsf g)$ accumulate at the thresholds $\nu_1=\pi^2/4$ and $\nu_2=\pi^2$ correspondingly; for details we refer to~\cite{Edward}. 
 
 Other examples of manifolds with axial analytic asymtotically cylindrical ends can be found in~\cite{KalvinII}, see also~\cite{KalvinI}.

\section{An approach to the complex scaling}\label{sCS}
This section presents a broad overview of our approach~\cite{KalvinII} to the complex scaling on manifolds with axial analytic asymptotically cylindrical ends. The approach originates from the one developed in~\cite{Hunziker} for $N$-body Schr\"{o}dinger operators.

We use the complex scaling $\Bbb R_+\ni x\mapsto
x+\lambda \mathsf s_R(x)$ along the axis  of the semi-cylinder $\Pi=\Bbb R_+\times\Omega$. Here $\mathsf s_R(x)=\mathsf s(x-R)$ is a scaling function with  a sufficiently large parameter $R>0$ and a smooth function $\mathsf s$  possessing the properties:
\begin{equation}\label{ab}
\begin{aligned}
&\mathsf s(x)=0 \text{ for all } x\leq 1,\\
&0\leq \mathsf s'(x)\leq 1 \text{ for all } x\in\Bbb R, \text{ and }
\mathsf s'(x)= 1 \text{ for large
 } x>0,
\end{aligned}
\end{equation}
where $\mathsf s'=\partial \mathsf s/\partial x$. The scaling parameter $\lambda$ takes its values in the disk
\begin{equation}\label{disk}
\mathcal D_\alpha=\{\lambda\in\mathbb C: |\lambda|<\sin\alpha< 1/\sqrt{2}\},
\end{equation}
where $\alpha<\pi/4$ is some angle for which the conditions of Definition~\ref{ACE} hold.
The function $\Bbb R_+\ni x\mapsto
x+\lambda \mathsf s_R(x)$ is invertible for all real $\lambda\in(-1,1)$, and
thus defines the selfdiffeomorphism
$$
\Pi\ni (x,\mathrm y)\mapsto\varkappa_\lambda (x,\mathrm y)=(x+\lambda \mathsf s_R(x),\mathrm y)\in \Pi,
$$
which scales the semi-cylinder $\Pi$ along its axis. We extend $\varkappa_\lambda$ to a selfdiffeomorphism of $\mathcal M$ by setting $\varkappa_\lambda(p)=p$ for all $p\in\mathcal M\setminus\Pi$. As a result we get
the Riemannian manifolds $(\mathcal M, \mathsf g_\lambda)$ parametrized by
 $\lambda\in(-1,1)$, where
$\mathsf g_\lambda=\varkappa_\lambda^*\mathsf g$ is the pullback
of the metric $\mathsf g$ by  $\varkappa_\lambda$.

\begin{figure}
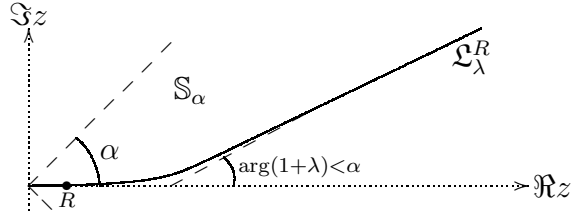
 \[ \xy0;/r.15pc/:
{\ar@{.>} (-30,0);(75,0)};(80,0)*{\Re z};
(-30,0);(0,30)**\dir{--};(-30,0);(-25,-5)**\dir{--};
(-22,0)*{\scriptstyle\bullet};(-22,-3)*{\scriptstyle R};
(-15,0);(-20,10)**\crv{(-16,7)};(-13,7)*{\alpha};(4,20)*{\Bbb
S_\alpha};
{\ar@{.>} (-30,-5);(-30,33)}; (-30,36)*{ \Im z};
(-30,0);(65,33)**\crv{(-5,0)&(5,4)&(15,9)&(36,19)};
(0,0);(28,15)**\dir{--};
(11,6);(13,0)**\crv{(14,3)};(27,4)*{\scriptstyle\arg(1+\lambda)<\alpha};
(63,27)*{ \mathfrak L^R_\lambda};
\endxy\]
\caption{The curve $\mathfrak
L^R_\lambda=\{z\in\Bbb C: z=x+\lambda \mathsf s_R(x),x\in\Bbb R_+
\}$ with $\lambda\in\mathcal D_\alpha$.}\label{fig+}
\end{figure}

Let  $\mathrm T'^*\Bbb S_\alpha$ be the holomorphic cotangent bundle $\{(z, c\,dz):z\in\Bbb S_\alpha, c\in\Bbb C\}$ of the sector $\Bbb S_\alpha=\{z\in\Bbb C:|\arg z|<\alpha<\pi/4\}$, where $dz=d\Re z+id\Im z$.
Consider the tensor field
\begin{equation}\label{TF}
 \mathfrak g_0 dz\otimes dz +2\mathfrak g_1\otimes dz+\mathfrak g_2\in C^\infty({\mathrm T}'^*\Bbb
S_\alpha\otimes{\mathrm T}^*\Omega )^{\otimes 2}
\end{equation}
with the analytic coefficients $\Bbb S_\alpha\ni z\mapsto \mathfrak g_k(z)\in C^\infty \mathrm T^*\Omega^{\otimes k}$, cf. Definition~\ref{ACE}. For all $\lambda\in\mathcal D_\alpha$ the complex scaling  defines the embedding
\begin{equation}\label{emb}
\mathrm T^*\Bbb R_+\ni\{x, a\,dx\}\mapsto \{x+\lambda
\mathsf s_R(x), a(1+\lambda \mathsf s'_R(x))^{-1} d z\}\in \mathrm T'^*\Bbb S_\alpha,
\end{equation}
where $|1+\lambda \mathsf s'_R(x)|>1-1/\sqrt{2}$, see Fig.~\ref{fig+}.
We identify the bundle $\mathrm T^*\Pi$ with the product $\mathrm T^*\Bbb R_+\otimes \mathrm T^*\Omega$ via the  product structure on $\Pi=\Bbb R_+\times\Omega$.
The embedding~\eqref{emb} together with~\eqref{TF} induces the tensor field
\begin{equation}\label{TFL}
\begin{aligned}
\mathsf g_\lambda\!\upharpoonright_\Pi=\mathfrak g^R_{0,\lambda}  dx\otimes dx +2\mathfrak g^R_{1,\lambda} \otimes dx +\mathfrak g^R_{2,\lambda}\in C^\infty {\mathrm T}^* \Pi^{\otimes 2},\\
\mathfrak g^R_{k,\lambda}(x)=(1+\lambda \mathsf s'_R(x))^{2-k}\mathfrak g_{k}(x+\lambda \mathsf s_R(x)),
\end{aligned}
\end{equation}
where
$\mathfrak g^R_{k,\lambda}(x)\in C^\infty \mathrm T^*\Omega^{\otimes k}$
are  smooth in $x\in\Bbb R_+$ and analytic in $\lambda\in\mathcal D_\alpha$ coefficients.
 Since $\supp
\mathsf s_R\cap (0,R)=\varnothing$, the equality $\mathsf g_\lambda\!\upharpoonright_{(0,R)\times\Omega}=\mathsf g\!\upharpoonright_{(0,R)\times\Omega}$ holds for all $\lambda\in\mathcal D_\alpha$, cf.~\eqref{split} and~\eqref{TFL}. Thanks to this we can smoothly extend $\mathsf g_\lambda\!\upharpoonright_{\Pi}$  to $\mathcal M$ by setting $\mathsf g_\lambda\!\upharpoonright_{\mathcal M\setminus\Pi}=\mathsf g\!\upharpoonright_{\mathcal M\setminus\Pi}$. As a result we obtain an analytic function
$$
\mathcal D_\alpha\ni\lambda\mapsto \mathsf g_\lambda \in C^\infty \mathrm T^*\mathcal M^{\otimes 2}.
$$

We consider the tensor field $\mathsf g_\lambda$  with  $\lambda\in\mathcal D_\alpha$ as a deformation of the metric $\mathsf g$ on $\mathcal M$ by means of the complex scaling.
Clearly, for $\lambda\in\mathcal D_\alpha\cap\Bbb R$ we have $\mathsf g_\lambda=\varkappa_\lambda^*\mathsf g$, and  $\mathsf g_0\equiv\mathsf g$.  By analyticity in $\lambda$ we conclude that $\mathsf g_\lambda$ is a symmetric tensor field. The Schwarz reflection principle gives $\overline{\mathsf g_\lambda}=\mathsf g_{\overline{\lambda}}$, where the bar stands for the complex conjugation. It must be stressed that the  tensor field $\mathsf g_\lambda$ with $\lambda\neq 0$ depends on  $R$, however we do not indicate this for brevity of notations.
As shown in~\cite{KalvinII}, the condition ii in Definition~\ref{ACE} 
implies that the tensor fields $\mathsf g_\lambda$ with $\lambda\in\mathcal D_\alpha$ are non-degenerate as $R>0$ is sufficiently large.

 It is well known that Riemannian metrics induce  musical isomorphisms between the tangent and cotangent bundles. Similarly, the non-degenerate tensor field $\mathsf g_\lambda$ induces a musical fiber isomorphism
$^\lambda\sharp:\Bbb C\mathrm T^*\mathcal M\to \Bbb C \mathrm T\mathcal  M$
between the complexified  bundles. For every $p\in\mathcal M$ the tensor field $\mathsf g_\lambda\in C^\infty\mathrm T^*\mathcal M^{\otimes 2}$  naturally defines a non-degenerate sesquilinear form $\mathsf g_\lambda^p[\cdot,\cdot]$ on $\Bbb C\mathrm T_p\mathcal M$.
 The isomorphism $^\lambda\sharp$ acts by the rule
$$
\Bbb C\mathrm T^*_p\mathcal M\ni\xi \mapsto {^\lambda\sharp} \xi\in \Bbb C \mathrm T_p\mathcal  M,
$$
where ${^\lambda\sharp} \xi$ is a unique vector satisfying the equality $ \xi\overline{\eta}=\mathsf g_\lambda^p [{^\lambda\sharp} \xi,\eta]$ for all $\eta \in \Bbb C \mathrm T_p\mathcal  M$. We extend the form $\mathsf g_\lambda^p[\cdot,\cdot]$ to the pairs $(\xi,\omega)\in {\Bbb  C{\mathrm T}^*_{p}\mathcal M}\times
{\Bbb C{\mathrm T}^*_{p}\mathcal M}$ by setting
\begin{equation*}
{\mathsf g_\lambda^p} [\xi,\omega]={\mathsf g_\lambda^p} [{^\lambda\sharp}\,\xi,{^{\overline{\lambda}}\sharp}\,\omega], \quad
\lambda\in\mathcal D_\alpha.
\end{equation*}

The function $\mathcal D_\alpha\ni\lambda\mapsto {\mathsf g_\lambda^p} [\xi,\omega]$ is analytic.
If $\lambda\in\mathcal D_\alpha$ is real, then ${\mathsf g_\lambda^p} [\cdot,\cdot]$ is the positive Hermitian form corresponding to the Riemannian metric $\mathsf g_\lambda$ on $\mathcal M$. In particular, ${\mathsf g_0^p} [\cdot,\cdot]\equiv \mathsf g^p [\cdot,\cdot]$. However, for a non-real $\lambda\in\mathcal D_\alpha$ we have $\overline{{\mathsf g_\lambda^p} [\xi,\omega]}={\mathsf g_{\overline \lambda}^p} [\omega,\xi]$, and the form ${\mathsf g_\lambda^p} [\cdot,\cdot]$ is not Hermitian. Nonetheless, on the differential one-forms the sesquilinear form ${\mathsf g_\lambda^p} [\cdot,\cdot]$ is sectorial and relatively bounded~\cite[Lemma 4.2]{KalvinII}. More precisely,
there exist some independent of  $p\in \mathcal M$ and $\lambda\in\mathcal D_\alpha$ angle $\vartheta< \pi/2$ and constant $\delta >0$, such that
\begin{equation}\label{new1}
 |\arg {\mathsf g_\lambda^p} [\xi,\xi]|\leq \vartheta, \quad \delta {\mathsf g^p} [\xi,\xi]\leq\Re{\mathsf g_\lambda^p} [\xi,\xi]\leq \delta^{-1} {\mathsf g^p} [\xi,\xi] \quad \forall \xi\in \Bbb C\mathrm T^*_p \mathcal M.
\end{equation}

By  $\operatorname{dvol}_{\lambda}$ with real $\lambda\in\mathcal D_\alpha$ we denote the volume form on the  Riemannian manifold  $(\mathcal M,\mathsf g_\lambda)$. The equality
$\varrho_\lambda\operatorname{dvol}_{\lambda}= \operatorname{dvol}_{0}$ defines a function  $\varrho_\lambda\in C^\infty(\mathcal M)$. Due to the properties of $\mathsf g_\lambda$ it turns out that  $\varrho_\lambda\in C^\infty(\mathcal M)$ is an analytic function of $\lambda\in\mathcal D_\alpha$ obeying the estimates
\begin{equation}\label{0}
0<c\leq|\varrho_\lambda(p)|\leq 1/c,\quad \mathsf g^p[d\varrho_\lambda,d\varrho_\lambda]<1/c,\quad p\in\mathcal M,
\end{equation}
where $c$ is independent of $p$ and $\lambda$.
 We introduce the deformed volume form
\begin{equation*}
\operatorname{dvol}_{\lambda}:=\frac 1 {\varrho_\lambda} \operatorname{dvol}_{0},\quad \lambda\in\mathcal D_\alpha,
\end{equation*}
and the  deformed global inner product
\begin{equation*}
(\xi,\omega)_\lambda=\int_{\mathcal M} \mathsf g_\lambda[\xi,\omega] \,\operatorname{dvol}_{\lambda}, \quad \xi,\omega\in
C_c^\infty\mathrm T^*\mathcal M^{\otimes k}.
\end{equation*}
Let us stress that  for non-real $\lambda\in\mathcal D_\alpha$ the deformed volume form is complex-valued, and the deformed inner product $(\xi,\omega)_\lambda=\overline{(\omega,\xi)_{\overline{\lambda}}}$ is not Hermitian.

Let $L^2 \mathrm T^*\mathcal M^{\otimes k}$ be the completion of the set $C_c^\infty \mathrm T^*\mathcal M^{\otimes k}$ with respect to the global inner product $(\cdot,\cdot)\equiv (\cdot,\cdot)_0$. The estimates~\eqref{new1}  together with the bounds on $\varrho_\lambda$ imply that the deformed global inner product $(\cdot,\cdot)_\lambda$ extends  to a uniformly bounded  non-degenerate form in $L^2 \mathrm T^*\mathcal M^{\otimes k}$ with $k=0,1$; i.e. for some independent of $\lambda\in\mathcal D_\alpha$ constant $c>0$  we have
\begin{equation}\label{3}
|(\xi,\omega)_\lambda|^2\leq c(\xi,\xi)(\omega,\omega),\quad \forall\xi,\omega\in L^2 \mathrm T^*\mathcal M^{\otimes k},k=0,1,
\end{equation}
and for any nonzero $\xi\in L^2 \mathrm T^*\mathcal M^{\otimes k}$  there exists  $\omega\in L^2 \mathrm T^*\mathcal M^{\otimes k}$, such that $(\xi,\omega)_\lambda\neq 0$.
If $\lambda\in\mathcal D_\alpha$ is real, then $(\cdot,\cdot)_\lambda$ coincides with the global inner product on the Riemannian manifold $(\mathcal M,\mathsf g_\lambda)$, and $\sqrt{(\cdot,\cdot)_\lambda}$ is an equivalent norm in $L^2\mathrm T^*\mathcal M^{\otimes k}$, $k=0,1$.

Let  ${^\lambda\!\Delta}:C_c^\infty(\mathcal M)\to C_c^\infty(\mathcal M)$ with $\lambda\in\Bbb R\cap\mathcal D_\alpha$ be the Laplacian on the Riemannian manifold $(\mathcal M,\mathsf g_\lambda)$. In the case $\partial\mathcal M\neq\varnothing$ we also consider the operator  $^\lambda\!{\partial}_\nu: C_c^\infty(\mathcal M)\to C_c^\infty(\partial\mathcal M)$ of the Neumann boundary condition on $(\mathcal M,\mathsf g_\lambda)$.  On a manifold with an axial analytic asymptotically cylindrical end the operators ${^\lambda\!\Delta}$ and  $^\lambda\!{\partial}_\nu$ extend by analyticity from $\Bbb R\cap\mathcal D_\alpha$  to all $\lambda\in\mathcal D_\alpha$, see~\cite{KalvinII}. Moreover, by taking a sufficiently large parameter $R>0$ we arrange the complex scaling so that the differential operator ${^\lambda\!\Delta}$ on $\mathcal M$ is strongly elliptic, and in the case $\partial\mathcal M\neq \varnothing$ the pair $\{{^\lambda\!\Delta}, {^\lambda\!{\partial}_\nu}\}$ obeys the Shapiro-Lopatinski\v{\i} condition on $\partial\mathcal M$; see \cite[Lemmas 6.2 and 6.3]{KalvinII}. Recall that for a strongly elliptic operator with the Dirichlet boundary
condition the Shapiro-Lopatinski\v{\i} condition is always fulfilled, e.g.~\cite{KozlovMazyaRossmann,Lions Magenes}.

Let $\{\mathscr U_j,\kappa_j\}$ be a finite atlas on $(\Omega,\mathfrak h)$, and let $y\in\Bbb R^n$  be a system of local coordinates in a neighborhood $\mathscr U_j$. If $\partial\Omega\cap\mathscr U_j\neq\varnothing$, then we assume in addition that  all $y$ in the image of the set $\partial\Omega\cap\mathscr U_j$ under the diffeomorphism $\kappa_j$ are of the form $y=(y',y_n)$ with $y'\in\Bbb R^{n-1}$ and $y_n\geq 0$, and the set $\partial\Omega\cap\mathscr U_j$ is defined by the equality $y_n=0$. In the coordinates $(x,y)$ on $\Pi$, where $x\in\Bbb R_+$ is the axial coordinate, the operators ${^\lambda\!\Delta}$ and  $^\lambda\!{\partial}_\nu$ with  $\lambda\in\mathcal D_\alpha$ admit the local  representations
\begin{equation}\label{LR}
 \begin{aligned}
 {^\lambda\!\Delta}&=-\frac {1} {\sqrt{\det \mathbf
 g_\lambda }}\nabla_{xy}\cdot \sqrt{\det \mathbf
 g_\lambda }\,\,\mathbf g^{-1}_\lambda
 \nabla_{xy},
 \\
{^\lambda\!\partial_\nu}&=\bigl(0,\dots,0,1/\sqrt{\mathbf
g^{-1}_{\lambda,nn}}\bigr)\mathbf
g^{-1}_\lambda \upharpoonright_{y_n=0}\nabla_{xy}\quad\text{ if }\quad  \mathscr U_j\cap\partial\Omega\neq \varnothing.
\end{aligned}
\end{equation}
Here  the complex symmetric matrix $\mathbf g_\lambda$ corresponds to the  representation of the tensor field $\mathsf g_\lambda$ in the coordinates $(x,y)$, and $\nabla_{xy}\equiv(\partial_x,\partial_{y_1},\dots,\partial_{y_n})^\intercal$. As shown in~\cite[Lemma 4.1]{KalvinII}, the inverse matrix $\mathbf g_\lambda^{-1}$ possesses the property
\begin{equation}\label{stab}
\sum_{|r|+q\leq 1}\|\partial^q_x\partial_y^r(\mathbf g^{-1}_\lambda(x,y)-\diag\{(1+\lambda)^{-2},\mathbf h^{-1}(y)\})\|_2\to 0\text{ as } x\to+\infty
\end{equation}
uniformly in $\lambda$ and $y$, where the matrix $\mathbf h(y)$ corresponds to the  representation of the metric $\mathfrak h$ in the local coordinates, and $\|\mathbf g\|_2=\sqrt{\sum_{\ell,m=0}^n |\mathbf g_{\ell m}|^2}$ is the matrix norm. Let us remark that  ${^\lambda\!\Delta}={^0\!\Delta}$ and  ${^\lambda\!{\partial}_\nu}={^0\!{\partial}_\nu}$ on $\mathcal M\setminus\Pi$. We also note that in the case of a general manifold $(\mathcal M,\mathsf g)$ with an asymptotically cylindrical end the representations~\eqref{LR} and the property~\eqref{stab} remain valid for $\lambda=0$ (and even for all real $\lambda\in(-1,1)$) as it  follows from Definition~\ref{ace}.

Consider ${^\lambda\!\Delta}$ as an unbounded operator in the Hilbert space $L^2(\mathcal M)$, initially defined on a dense in $L^2(\mathcal M)$ core $\mathbf C({^\lambda\!\Delta})$.
\begin{defn}\label{d1} In the case
$\partial\mathcal M=\varnothing$  we take  $\mathbf C({^\lambda\!\Delta})\equiv C^\infty_c(\mathcal M)$. In the case of the Neumann (resp. Dirichlet) Laplacian $\Delta$ the core $\mathbf C({^\lambda\!\Delta})$ consists of the functions $u \in
C^\infty_c({\mathcal M})$ satisfying the deformed Neumann boundary condition  ${^\lambda\!{\partial}_\nu} u=0$ (resp. the Dirichlet boundary condition
$u\!\upharpoonright_{\partial \mathcal M}=0$).
\end{defn}
 The operator ${^\lambda\!\Delta}$ with $\lambda\in\mathcal D_\alpha$ is a deformation of the Laplacian $\Delta\equiv{^0\!\Delta}$ by means of the complex scaling.
In general, the operator ${^\lambda\!{\partial}_\nu}$ and the core $\mathbf C({^\lambda\!\Delta})$ depend on the scaling parameter $\lambda\in\mathcal D_\alpha$.

Let $d: C_c^\infty(\mathcal M)\to C_c^\infty \mathrm T^*\mathcal M$ be the exterior derivative. We introduce the sesquilinear quadratic form
\begin{equation*}
\mathsf q_\lambda[u,v]=\bigl(du,d(\overline{\varrho_\lambda}v)\bigr)_\lambda,\quad u,v\in\mathbf C(\mathsf q),\quad\lambda\in\mathcal D_\alpha,
\end{equation*}
on a core $\mathbf C(\mathsf q)$.
\begin{defn}\label{form core} In the case  $\partial\mathcal M=\varnothing$, and also in the case of the Neumann Laplacian,  we take $\mathbf C(\mathsf q)\equiv C_c^\infty(\mathcal M)$. In the case of the Dirichlet Laplacian the core  $\mathbf C(\mathsf q)$ consists of the functions $u\in C_c^\infty (\mathcal M)$ with $u\!\upharpoonright_{\partial\mathcal M}=0$.
\end{defn}
As it follows from the Green identity
\begin{equation}\label{1}
(du,d v)_\lambda=({^{\lambda}\!\Delta u}, v)_\lambda,\quad u\in \mathbf C({^\lambda\!\Delta}), v\in\mathbf C(\mathsf q),
\end{equation}
to the unbounded operator ${^\lambda\!\Delta}$ in the Hilbert space $L^2(\mathcal M)$ there corresponds the quadratic form  $\mathsf q_\lambda[\cdot,\cdot]$.
Clearly, $\mathsf q_0[\cdot,\cdot]=(d\cdot,d\cdot)$  is the nonnegative quadratic form of the Laplacian $\Delta$. Below we formulate a result from~\cite{KalvinII}, for the proof we refer to~\cite[Proposition 5.3]{KalvinII}.
\begin{prop}\label{p1} Assume that $(\mathcal M,\mathsf g)$ is a manifold with an axial analytic asymptotically cylindrical end. Let $\mathbf D(\mathsf q)$ be the Hilbert space, introduced as the completion of the core $\mathbf C(\mathsf q)$ with respect to the norm $\sqrt{(du,du)+\|u\|^2}$. Then the following assertions hold.  \begin{itemize}
\item[i.] The  unbounded quadratic form  $\mathsf q_\lambda [\cdot,\cdot]$ with the domain $\mathbf D(\mathsf q)$ is densely defined and closed for all $\lambda\in\mathcal D_\alpha$.
\item[ii.] The form $\mathsf q_\lambda[\cdot,\cdot]$ is sectorial and relatively bounded. More precisely, for all $\lambda\in\mathcal D_\alpha$ and    $u\in\mathbf D(\mathsf q)$ the estimates
 \begin{equation*}
 |\arg(\mathsf q_\lambda [u,u]+a \|u\|^2)|\leq \vartheta<\pi/2,
\end{equation*}
\begin{equation*}
b (du,du)\leq \Re \mathsf q_\lambda [u,u]+a\|u\|^2,\quad
\Re \mathsf q_\lambda [u,u]\leq b^{-1} \bigl((du,du)+\|u\|^2\bigr)
\end{equation*}
hold with some constants $\vartheta$ and $a,b>0$, which are independent of $\lambda$ and $u$.
\item[iii.] For any $u\in \mathbf D(\mathsf q)$ the function $\mathcal D_\alpha\ni\lambda\mapsto \mathsf q_\lambda[u,u]$  is analytic.
\end{itemize}
\end{prop}

\section{Conjugated operator and its essential spectrum}\label{s5}

In order to study exponential decay of eigenfunctions we consider the operator ${^\lambda\!\Delta}$ conjugated with an exponent and use the dilation analytic techniques, similar to those we meet in the theory of $N$-body Schr\"{o}dinger operators, e.g.~\cite[Chapter XIII.11]{Simon Reed iv}.

Let  $\mathsf s$ be a smooth function on the semi-cylinder $\Pi$, which depends only on the axial variable $x\in\Bbb R_+$ and possesses  the properties~\eqref{ab}. We extend $\mathsf s$ to a smooth function on $\mathcal M$ by setting $\mathsf s\!\upharpoonright_{\mathcal M\setminus\Pi}\equiv 0$. Consider the conjugated operator ${^\lambda\!\Delta}_\beta=e^{-\beta\mathsf s}\,{^\lambda\!\Delta}\,e^{\beta\mathsf  s}$ with the  parameter $\beta\in\Bbb C$, where  $e^{\beta\mathsf  s}$ is the operator of multiplication by the exponent. The unbounded operator  ${^\lambda\!\Delta}_\beta$ in $L^2(\mathcal M)$ is initially defined on the dense in $L^2(\mathcal M)$ core
\begin{equation}\label{4}
\mathbf C({^\lambda\!\Delta}_\beta)=\{ u: e^{\beta\mathsf s} u\in \mathbf C({^\lambda\!\Delta})\};
\end{equation}
here $\mathbf C({^\lambda\!\Delta})$ is the same as in Definition~\ref{d1}. The core $\mathbf C({^\lambda\!\Delta}_\beta)$  depends on the  parameters $\beta\in\Bbb C$ and $\lambda\in\mathcal D_\alpha$. From~\eqref{1} and~\eqref{4} we get
$$
(e^{-\beta\mathsf s}\,{^\lambda\!\Delta}\,e^{\beta \mathsf s}u,v)=\bigl( d(e^{\beta \mathsf s} u),d(\overline{e^{-\beta \mathsf s}\varrho_\lambda}v)\bigr)_\lambda,\quad u\in\mathbf C({^\lambda\!\Delta}_\beta ), v\in\mathbf C(\mathsf q).
$$
Thus to the operator ${^\lambda\!\Delta}_\beta$ there corresponds the quadratic form
$$
\mathsf q_\lambda^\beta[u,v]=\bigl( d(e^{\beta \mathsf s} u),d(\overline{e^{-\beta \mathsf s}\varrho_\lambda}v)\bigr)_\lambda, \quad u,v\in \mathbf C(\mathsf q).
$$

\begin{lem}~\label{l1} Assume that $(\mathcal M,\mathsf g)$ is a manifold with an axial analytic asymptotically cylindrical end. Then the difference $\mathsf q^\beta_\lambda[\cdot,\cdot]- \mathsf q_\lambda[\cdot,\cdot]$ has an arbitrarily small uniform in $\lambda\in\mathcal D_\alpha$ relative bound  with respect to the form $\mathsf q_\lambda[\cdot,\cdot]$.  More precisely, for all $\lambda\in\mathcal D_\alpha$ and $u\in\mathbf C(\mathsf q)$ the estimate
$$
|\mathsf q_\lambda^\beta[u,u]-\mathsf q_\lambda[u,u]|\leq \varepsilon |\mathsf q_\lambda[u,u]|+ C(|\beta|,\varepsilon)\|u\|^2
$$
is valid, where $\varepsilon>0$ is arbitrarily small,  and the constant $C(|\beta|,\varepsilon)$  depends on $|\beta|$ and $\varepsilon$, but not on $u$ or $\lambda$.

If $(\mathcal M,\mathsf g)$ is a general manifold with an asymptotically cylindrical end, then the assertion remains valid for $\lambda=0$.
\end{lem}

\begin{pf} Let $u\in\mathbf C(\mathsf q_0)$.
After simple calculations  we obtain the equality
\begin{equation}\label{01}
\begin{aligned}
\mathsf q_\lambda^\beta[u,u]-&\mathsf q_\lambda[u,u]=\beta\bigl(u\, d\mathsf s,d(\overline{\,\varrho_\lambda}u)\bigr)_\lambda\\&
-\beta\bigl(\varrho_\lambda\, du,u\,d\mathsf s\bigr)_\lambda
-\beta^2\bigl(\varrho_\lambda u\,d\mathsf s,u\,d\mathsf s\bigr)_\lambda.
\end{aligned}
\end{equation}
We will rely on~\eqref{3}. Then for an arbitrarily small $\epsilon>0$ the first term in the right hand side of~\eqref{01} meets the estimates
$$
\bigl|\bigl(u\, d\mathsf s,d(\overline{\,\varrho_\lambda}u)\bigr)_\lambda\bigr|\leq c \bigl(u\,d\mathsf s,u\,d\mathsf s\bigr)^{1/2}\bigl(d(\overline{\varrho_\lambda}u),d(\overline{\varrho_\lambda}u)\bigr)^{1/2}
$$
$$\leq  c \epsilon \bigl(d(\overline{\varrho_\lambda}u),d(\overline{\varrho_\lambda}u)\bigr)+c\epsilon^{-1}\bigl(u\,d\mathsf s,u\,d\mathsf s\bigr)
$$
$$
\leq  2c\epsilon(u\,d\overline{\varrho_\lambda},u\,d\overline{\varrho_\lambda})+2c\epsilon(|\varrho_\lambda|^2du,du)+c\epsilon^{-1}\bigl(u\,d\mathsf s,u\,d\mathsf s\bigr).
$$
For the summands in the last line we have
$$
(u\,d\overline{\varrho_\lambda},u\,d\overline{\varrho_\lambda})\leq \sup_{p\in\mathcal M}\mathsf g^p[d{\varrho_\lambda},d{\varrho_\lambda}]\|u\|^2,
\quad(|\varrho_\lambda|^2du,du)\leq\sup_{p\in\mathcal M}|\varrho_\lambda(p)|^2(du,du),
$$
$$
\bigl(u\,d\mathsf s,u\,d\mathsf s\bigr)\leq \sup_{p\in\mathcal M}\mathsf g^p[d\mathsf s,d\mathsf s] \|u\|^2.
$$
The second term in the right hand side of~\eqref{01} obeys the estimate
$$
\bigl|\bigl(\varrho_\lambda\, du,u\,d\mathsf s\bigr)_\lambda\bigr|\leq c\epsilon\bigl(|\varrho_\lambda|^2\, du, du\bigr)+c\epsilon^{-1}\bigl(u\,d\mathsf s,u\,d\mathsf s\bigr).
$$
Finally, for the last term in~\eqref{01} we get
$$
\bigl|\bigl(\varrho_\lambda u\,d\mathsf s,u\,d\mathsf s\bigr)_\lambda\bigr|\leq c\bigl(\varrho_\lambda u\,d\mathsf s,\varrho_\lambda u\,d\mathsf s\bigr)^{1/2}\bigl(u\,d\mathsf s,u\,d\mathsf s\bigr)^{1/2}
$$
$$
\leq c\bigl(\sup_{p\in\mathcal M} |\varrho_\lambda(p)|^2\bigr)^{1/2} \sup_{p\in\mathcal M}\mathsf g^p[d\mathsf s,d\mathsf s]\|u\|^2.
$$

Observe that  $d\mathsf s\!\upharpoonright_{\mathcal M\setminus\Pi}=0$ and $d\mathsf s\!\upharpoonright_{\Pi}=\mathsf s'\, dx$, where $\mathsf s'(p)=\mathsf s'(x)\leq 1$. Due to stabilization of the metric $\mathsf g$ to the product metric $dx\otimes dx+\mathfrak h$ at infinity,  we have $\sup_{p\in\mathcal M}\mathsf g^p[d\mathsf s,d\mathsf s]\leq C$. Now the assertion follows from  the obtained estimates combined with~\eqref{0},~\eqref{01}, and Proposition~\ref{p1}.ii.
\qed
\end{pf}

\begin{prop}\label{q beta} Assume that $(\mathcal M,\mathsf g)$ is a manifold with an axial analytic asymptotically cylindrical end. Let $\lambda\in\mathcal D_\alpha$ and $\beta\in\Bbb C$. Then the following assertions hold.  \begin{itemize}
\item[i.] The  unbounded sesquilinear form  $\mathsf q_\lambda^\beta [\cdot,\cdot]$ with the domain $\mathbf D(\mathsf q)$ is densely defined and closed.
\item[ii.] The form $\mathsf q^\beta_\lambda[\cdot,\cdot]$ is sectorial and relatively bounded. More precisely, for all $u\in\mathbf D(\mathsf q)$ the estimates
 \begin{equation*}
 |\arg(\mathsf q^\beta_\lambda [u,u]+a \|u\|^2)|\leq \vartheta,
\end{equation*}
\begin{equation*}
b (du,du)\leq \Re \mathsf q^\beta_\lambda [u,u]+a\|u\|^2,\quad
\Re \mathsf q^\beta_\lambda [u,u]\leq b^{-1} \bigl((du,du)+\|u\|^2\bigr)
\end{equation*}
hold with some angle $\vartheta<\pi/2$ and some positive constants $a$ and $b$, which may depend on $\beta$ and $\vartheta$, but not on $\lambda$ or $u$.
\item[iii.] For any $u\in \mathbf D(\mathsf q)$ and $\lambda\in\mathcal D_\alpha$ the function $\Bbb C\ni\beta\mapsto \mathsf q^\beta_\lambda[u,u]$  is analytic.
\end{itemize}

In the case of a general manifold with an asymptotically cylindrical end the assertions remain valid for $\lambda=0$.
\end{prop}
\begin{pf} The first two assertions are direct  consequences of Lemma~\ref{l1} and Proposition~\ref{p1}; all necessary basic facts from the theory of sesquilinear forms can be found  e.g. in~\cite[Chapter VI]{Kato}. The last assertion is an immediate consequence of the equality~\eqref{01}.
\qed\end{pf}

As is known~\cite[Chapter VI.2.1]{Kato}, there is a one-to-one correspondence between the set of all densely defined closed sectorial sesquilinear forms and the set of all m-sectorial operators.
(Here and elsewhere m-sectorial means that the numerical range $\{({\mathsf A}u,u): u\in\mathbf D(\mathsf A) \}$
and the spectrum of a  closed operator $\mathsf A$ with domain $\mathbf D(\mathsf A)$ lie in the sector $\{z\in\Bbb C:|\arg (z+a)|\leq \vartheta\}$ with some $\vartheta<\pi/2$ and $a>0$.) Thus Proposition~\ref{q beta} implies that the Friedrichs extension of the operator ${^\lambda\!\Delta}_\beta$, initially defined on the dense in $L^2(\mathcal M)$ core $\mathbf C({^\lambda\!\Delta}_\beta)$, is m-sectorial.   In particular, the Friedrichs extension of the Laplacian $\Delta\equiv {^0\!\Delta}_0$ is a nonnegative selfadjoint operator.

\begin{figure}[h]
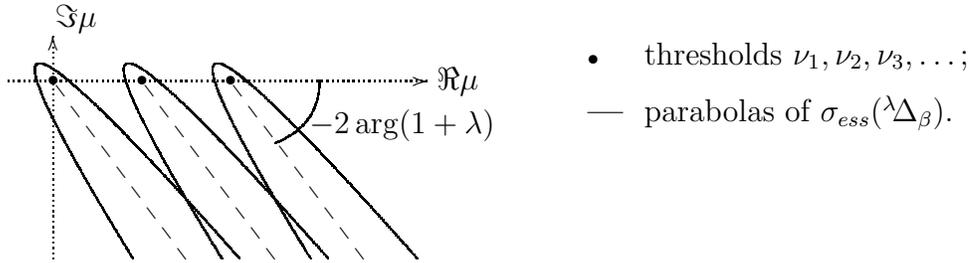

\[
\xy (0,0)*{\xy0;/r.28pc/:{\ar@{.>}(0,10);(50,10)*{\ \Re \mu}};
{\ar@{.>}(5,-10);(5,15)}; (7,17)*{\ \Im \mu};
{\ar@{--}(5,10)*{\scriptstyle\bullet};(20,-10)};(14,-10);(26,-10)**\crv{(-13,34)};
{\ar@{--}(15,10)*{\scriptstyle\bullet};(30,-10)};(24,-10);(36,-10)**\crv{(-3,34)};
{(25,10)*{\scriptstyle\bullet};(40,-10)*{} **\dir{--}};
(35,10);(30,3)**\crv{(35,5)*{\quad\quad\quad\quad\quad\
-2\arg(1+\lambda)}};(34,-10);(46,-10)**\crv{(7,34)};
\endxy};
(70,0)*{\xy (20,20)*{
\begin{array}{ll}
     {\scriptstyle\bullet } & \text{ thresholds $\nu_1,\nu_2,\nu_3,\dots$;} \\
    \text{\bf---} & \text{ parabolas of $\sigma_{ess}({^\lambda\!\Delta}_\beta)$.}
\end{array}};
\endxy};
\endxy
\]
\caption{Essential spectrum of the operator
${^\lambda\!\Delta}_\beta$ for $\Im\lambda>0$ and $\beta\gtrless 0$.}\label{fig5}
\end{figure}

Consider the domain $\mathbf D({^\lambda\!\Delta}_\beta)$  of the  m-sectorial operator ${^\lambda\!\Delta}_\beta$   as a Hilbert
space introduced as the completion of the core~\eqref{4} with respect to the graph norm
$\sqrt{\|\cdot\|^2+\|{^\lambda\!\Delta}_\beta\cdot\|^2}$.
We say that $\mu$ is a point of the essential spectrum
$\sigma_{ess}({^\lambda\!\Delta}_\beta)$, if the bounded operator
\begin{equation}\label{5}
{^\lambda\!\Delta}_\beta-\mu :\mathbf D({^\lambda\!\Delta}_\beta)\to L^2(\mathcal M)
\end{equation}
is not  Fredholm. (Recall that  a bounded linear operator is said to be  Fredholm, if  its kernel  and cokernel
 are finite-dimensional, and the
range  is closed.) In the next proposition we localize the essential spectrum of the operator ${^\lambda\!\Delta}_\beta$. For $\beta=0$ the proposition was already proven in~\cite[Theorem 6.1]{KalvinII}, the general case $\beta\in\Bbb C$ is very similar. In the proof we use methods of the theory of non-homogeneous elliptic boundary value problems~\cite{Lions Magenes,MP2}, see also~\cite{KozlovMaz`ya,KozlovMazyaRossmann}.
\begin{prop}\label{ess}  Assume that $(\mathcal M,\mathsf g)$ is a manifold with an axial analytic asymptotically cylindrical end.
Let $\lambda\in \mathcal D_\alpha$ and $\beta\in\Bbb C$.  Then $\mu\in \sigma_{ess}({^\lambda\!\Delta}_\beta)$, if and only if the equality
\begin{equation}\label{eq9}
\nu_j-\mu=(1+\lambda)^{-2}(\beta+i\xi)^2
\end{equation}
holds for some $j\in\mathbb N$  and some $\xi\in\Bbb R$,
where $\{\nu_j\}_{j=1}^\infty$ is the set of thresholds of the Laplacian $\Delta$ on $(\mathcal M,\mathsf g)$.

The spectrum $\sigma_{ess}({^\lambda\!\Delta}_\beta)$ is depicted on Fig.~\ref{fig5}. In the case $\beta=0$ the parabolas of the essential spectrum collapse to the dashed rays originating from every threshold $\nu_j$, and we obtain the essential spectrum $\sigma_{ess}({^\lambda\!\Delta})$ of ${^\lambda\!\Delta}\equiv{^\lambda\!\Delta}_0$.

In the case of a general manifold with an asymptotically cylindrical end the assertion remains valid for $\lambda=0$.
\end{prop}

\begin{pf} As it was already mentioned, by taking a sufficiently large $R$ in the function $\mathsf s_R(x)=\mathsf s(x-R)$, we arrange the complex scaling so that the deformation ${^\lambda\!\Delta}$ is a strongly elliptic operator on $\mathcal M$, and the pair $\{{^\lambda\!\Delta}, {^\lambda\!{\partial}_\nu}\}$ meets the Shapiro-Lopatinski\v{\i} condition on $\partial\mathcal M$, if $\partial\mathcal M\neq\varnothing$. The principal symbols of the operators ${^\lambda\!\Delta}_\beta$ and ${^\lambda\!\Delta}$  are coincident, as well as the principal symbols  of $e^{-\beta\mathsf s}\,{^\lambda\!{\partial}_\nu \,e^{\beta\mathsf  s}}$ and ${^\lambda\!{\partial}_\nu}$. Hence the differential operator ${^\lambda\!\Delta}_\beta$ is strongly elliptic on $\mathcal M$, and in the case $\partial\mathcal M\neq\varnothing$ the pair $\{{^\lambda\!\Delta}_\beta, e^{-\beta\mathsf s}\,{^\lambda\!{\partial}_\nu \,e^{\beta\mathsf  s}}\}$ meets the Shapiro-Lopatinski\v{\i} condition on $\partial\mathcal M$. All other difficulties related to the appearance of the deformed operator of the Neumann boundary conditions $e^{-\beta\mathsf s}\,{^\lambda\!{\partial}_\nu \,e^{\beta\mathsf  s}}$ on $\partial \mathcal M$ can be handled exactly in the  same way as in~\cite[Theorem 6.1]{KalvinII}. For this reason here we consider only the case of the Dirichlet Laplacian on a manifold $(\mathcal M,\mathsf g)$ with non-compact boundary $\partial \mathcal M$. The case $\partial\mathcal M=\varnothing$ is similar. As in~\cite{KalvinI,KalvinII} we will rely on the  Peetre's lemma:
\begin{itemize}
\item[]{\it Let $\mathcal X,\mathcal Y$ and $\mathcal Z$ be Banach spaces, where $\mathcal X$ is compactly embedded into $\mathcal Z$. Furthermore, let $\mathcal L$ be a linear continuous operator from $\mathcal X$ to $\mathcal Y$. Then the next two assertions are equivalent: (i) the range of $\mathcal L$ is closed in $\mathcal Y$ and $\dim \ker \mathcal L<\infty$, (ii) there exists a constant $C$, such that
    \begin{equation}\label{coercive}
    \|u\|_{\mathcal X}\leq C(\|\mathcal L u\|_{\mathcal Y}+\|u\|_{\mathcal Z})\quad \forall u\in \mathcal X.
    \end{equation}}
\end{itemize}
For the proof of this lemma we refer to~\cite{Peetre},~\cite[Lemma~5.1]{Lions Magenes}.

{\it Sufficiency.}
Here we assume that the spectral parameter $\mu$ does not meet the condition~\eqref{eq9}, and establish an estimate of type~\eqref{coercive} for the operator~\eqref{5}.

Introduce the Sobolev space ${\overset{\circ}H}\vphantom{H}^\ell(\Bbb R\times\Omega)$ of functions on the infinite cylinder $\Bbb R\times \Omega$ as the completion of the set $C_0^\infty(\Bbb R\times\Omega)$ with respect to the norm
\begin{equation}\label{norm}
\|u\|_{H^\ell(\Bbb R\times\Omega)}=\Bigl (\int_{\Bbb R}\sum_{r\leq \ell}\|\partial^r_x u\|^2_{ H^{\ell-r}(\Omega)}\,dx\Bigr)^{1/2},
\end{equation}
where $H^{\ell-r}(\Omega)$ is the Sobolev space of functions on the compact manifold $\Omega$.
Let  $L^2(\Bbb R\times\Omega)$ be the space induced by the product metric $dx\otimes dx+\mathfrak h$ on $\Bbb R\times\Omega$. By applying the Fourier transform $\mathcal F_{x\mapsto \xi}$ we pass from the continuous operator
\begin{equation}\label{opo}
\Delta_\Omega-(1+\lambda)^{-2}(\partial_x+\beta)^2 -\mu: {\overset{\circ}H}\vphantom{H}^2(\Bbb R\times\Omega)\to L^2(\Bbb R\times\Omega)
\end{equation}
of the Dirichlet boundary value problem in the infinite cylinder $\Bbb R\times\Omega$ to the Dirichlet Laplacian  $\Delta_\Omega+(1+\lambda)^{-2}(\beta+i\xi)^2-\mu$ with the spectral parameter $\mu-(1+\lambda)^{-2}(\beta+i\xi)^2$. Assume that $\mu$ does not satisfy the equality~\eqref{eq9} for any $j\in\Bbb N$ and $\xi\in \Bbb R$ or, equivalently, assume that  for any $\xi\in\Bbb R$ the number $\mu-(1+\lambda)^{-2}(\beta+i\xi)^2$  is not an eigenvalue $\nu_j$ of the Dirichlet Laplacian $\Delta_\Omega$ on $(\Omega,\mathfrak h)$.  Then a known argument, see e.g.~\cite[Theorem 4.1]{MP2} or \cite[Theorem 5.2.2]{KozlovMazyaRossmann} or \cite[Theorem~2.4.1]{KozlovMaz`ya}, shows that the operator~\eqref{opo} realizes an isomorphism. In particular, the estimate
\begin{equation}\label{op}
\|u\|_{ {H}^2(\Bbb R\times\Omega)}\leq C\|(\Delta_\Omega-(1+\lambda)^{-2}(\partial_x+\beta)^2 -\mu)u\|_{L^2(\Bbb R\times\Omega)}
\end{equation}
is valid with an independent of $u\in {\overset{\circ}H}\vphantom{H}^2(\Bbb R\times\Omega)$ constant $C=C(\mu,\lambda,\beta)>0$.

 Let $\chi_T(x)=\chi(x-T)$, where $\chi\in  C^\infty(\Bbb R)$ is a cutoff function, such that
$\chi(x)=1$ for $x\geq -3$ and $\chi(x)=0$ for $x\leq -4$. As a consequence of the stabilization of the tensor field $\mathsf g_\lambda$ to $(1+\lambda)^{2}dx\otimes dx+\mathfrak h$ at infinity, the operator ${^\lambda\!\Delta}_\beta$ stabilizes to $\Delta_\Omega-(1  +\lambda)^{-2}(\partial_x+\beta)^2$ at infinity in the sense that
$$
\|({^\lambda\!\Delta}_\beta-\Delta_\Omega+(1  +\lambda)^{-2}(\partial_x+\beta)^2)\chi_T u\|_{ L^2(\Bbb R\times\Omega)}\leq c(T)\|\chi_T u\|_{{H}^2(\Bbb R\times\Omega)},
$$
 where $c(T)\to 0$ as $T\to+\infty$, cf.~\eqref{LR} and~\eqref{stab}.   This together with~\eqref{op} implies that for a sufficiently large fixed $T=T(\mu,\lambda,\beta)>0$ the estimate
\begin{equation}\label{einf}
 \|\chi_Tu\|_{{H}^2(\Bbb R\times\Omega)} \leq \texttt{C}\|({^\lambda\!\Delta}_\beta -  \mu) \chi_T u\|_{ L^2(\Bbb R\times\Omega)}
\end{equation}
  holds, where the constant $\texttt{C}=(1/C- c(T))^{-1}>0$ may depend on $\mu$, $\lambda$, and $\beta$, but not on $u\in {\overset{\circ}H}\vphantom{H}^2(\Bbb R\times\Omega)$.

Without loss of generality we can  assume that $(0,T)\times\Omega\subset\mathcal M_c$, cf.~Fig~\ref{manif}. If it is not the case, then we take a larger smooth compact manifold $\mathcal M_c$, inserting the cylinder $(0,T)\times\Omega$ instead of the part $(0,1)\times\Omega$ of $\mathcal M_c$; recall that $(0,1)\times\Omega\subset\mathcal M_c\cap\Pi$ by our assumptions.

Let $\rho,\varsigma\in C_c^\infty(\mathcal M)$ be some cutoff functions, such that $\rho= 1$ on $\mathcal M\setminus(T-2,\infty)$ and $\rho= 0$ on $(T-1,\infty)\times\Omega$, while  $\varsigma\rho=\rho$ and $\supp\varsigma\subset\mathcal M_c$. As the operator ${^\lambda\!\Delta}_\beta$ is a strongly elliptic operator on $\mathcal M$,  the local coercive estimate
\begin{equation}\label{ecomp}
\|\rho u\|_{H^2(\mathcal M_c)}\leq C(\|\varsigma {^\lambda\!\Delta}_\beta u\|+\|\varsigma u\|)
\end{equation}
holds for all  $u\in C_0^\infty(\mathcal M)$. We write the estimate~\eqref{einf} for $u\in C_0^\infty(\mathcal M)$ in the form
$$
\|\chi_Tu\|_{{H}^2(\Bbb R\times\Omega)} \leq \texttt{C} \bigl(\|\chi_T(^\lambda\!\Delta_\beta -  \mu)u\|_{L^2(\Bbb R\times\Omega)}+\|[{^\lambda\!\Delta}_\beta,\chi_T]u\|_{ L^2(\Bbb R\times\Omega)}
\bigr),
$$
where the commutator $[{^\lambda\!\Delta}_\beta,\chi_T]$ is  equal to zero outside of the set $(T-5,T-2)\times\Omega\subset\mathcal M_c$. Since $\rho=1$ on this set, we get the estimate
$$
\|[{^\lambda\!\Delta}_\beta,\chi_T]u\|_{ L^2(\Bbb R\times\Omega)}\leq C \|\rho u\|_{H^2(\mathcal M_c)}.
$$
Due to stabilization of $\mathsf g$ at infinity to the product metric  $dx\otimes dx+\mathfrak h$  we  have
$$
\|\chi_T F\|_{L^2(\Bbb R\times\Omega)}^2=\int_{\Bbb R_+}\|\chi_T F\|^2_{L^2(\Omega)}\,dx\leq C\|F\|^2\quad \forall F\in L^2(\mathcal M).
$$

Introduce the Sobolev space ${\overset{\circ}H}\vphantom{H}^2(\mathcal M)$ as the completion of the set $C_0^\infty(\mathcal M)$ with respect to the norm
$$
\| u\|_{{H}^2(\mathcal M)}:=\|\chi_Tu\|_{ H^{2}(\Bbb R\times\Omega)}+\|\rho u\|_{H^2(\mathcal M_c)}.
$$
Then from the last four estimates it follows that
\begin{equation}\label{coercive!}
\| u\|_{{H}^2(\mathcal M)}\leq C\bigl(\|({^\lambda\!\Delta}_\beta -  \mu)u\|+\|\varsigma u\|\bigr),
\end{equation}
where the constant $C$ depends on $\lambda$ and $\mu$, but not on $u\in{\overset{\circ}H}\vphantom{H}^2(\mathcal M)$.
 We also have the estimate
$$
 \|{^\lambda\!\Delta}_\beta u\|+\|u\|\leq c \| u\|_{{H}^2(\mathcal M)}\quad\forall u\in {\overset{\circ}H}\vphantom{H}^2(\mathcal M).
$$
 This together with the estimate~\eqref{coercive!} implies that the spaces ${\overset{\circ}H}\vphantom{H}^2(\mathcal M)$ and $\mathbf D({^\lambda\!\Delta}_\beta)$ are coincident and their norms are equivalent.

Let $\mathsf w$ be a bounded rapidly decreasing at infinity positive function on $\mathcal M$, such that the embedding of ${\overset{\circ}H}\vphantom{H}^2(\mathcal M)$ into the weighted space $L^2(\mathcal M,\mathsf w)$ with the norm $\|\mathsf w \cdot\|$ is compact. As a consequence of~\eqref{coercive!} we obtain the estimate
\begin{equation}\label{coer}
\| u\|_{{H}^2(\mathcal M)}\leq {\mathrm C}\bigl(\|({^\lambda\!\Delta}_\beta -  \mu)u\|+\|\mathsf w u\|\bigr)\quad\forall u\in{\overset{\circ}H}\vphantom{H}^2(\mathcal M)
\end{equation}
of type~\eqref{coercive}. Then by the Peetre's lemma the range of the continuous operator ${^\lambda\!\Delta}_\beta -  \mu:{\overset{\circ}H}\vphantom{H}^2(\mathcal M)\to L^2(\mathcal M)$ is closed and the kernel is finite dimensional.
In order to see that the cokernel of this operator is finite dimensional, one can apply a similar argument to the  adjoint operator $^\lambda\!\Delta^*_\beta=e^{\overline{\beta\mathsf s}}\frac {1} {\varrho_{\overline{\lambda}}}{^{\overline{\lambda}}\!\Delta}{\varrho_{\overline{\lambda}}}e^{\overline{-\beta\mathsf s}}$. The operator $^\lambda\!\Delta^*_\beta$ stabilizes at infinity to the operator $\Delta_\Omega-(1+\overline{\lambda})^{-2}(\overline{\beta}-\partial_x)^2$. If $\mu$ does not meet the condition~\eqref{eq9}, this allows to deduce the estimate
$$
\| u\|_{{H}^2(\mathcal M)}\leq \mathrm C\bigl(\|({^\lambda\!\Delta}_\beta-\mu)^*u\|+\|\mathsf w u\|\bigr)\quad\forall u\in{\overset{\circ}H}\vphantom{H}^2(\mathcal M),
$$
which implies that the cokernel  $\operatorname{coker}({^\lambda\!\Delta} -  \mu)=\ker (^\lambda\!\Delta-\mu)^*$ is finite dimensional. Thus the deformed Dirichlet Laplacian ${^\lambda\!\Delta}_\beta-\mu:\mathbf D({^\lambda\!\Delta}_\beta)\to L^2(\mathcal M)$ is Fredholm, if $\mu$ does not meet the condition~\eqref{eq9}.

{\it Necessity.} Now we assume that $\mu$ meets the condition~\eqref{eq9} for some $j$, and show that 
 the operator ${^\lambda\!\Delta}_\beta -  \mu:{\overset{\circ}H}\vphantom{H}^2(\mathcal M)\to L^2(\mathcal M)$ is not Fredholm.  By the Peetre's lemma it suffices to find a sequence $\{u_\ell\}_{\ell=1}^\infty$ of functions $u_\ell\in{\overset{\circ}H}\vphantom{H}^2(\mathcal M)$ violating the estimate~\eqref{coer}.

Let $\chi$ be a smooth cutoff
function on the real line, such that $\chi(x)=1$ for $|x-3|\leq 1$,
and $\chi(x)=0$ for $|x-3|\geq 2$. Consider the functions
\begin{equation}\label{test}
u_\ell(x,\mathrm
y)=\chi(x/\ell)\exp\bigl({i(1+\lambda)x\sqrt{\mu-\nu_j}}-\beta x\bigr)\Phi(\mathrm
y),\quad (x,\mathrm y)\in \Bbb R\times\Omega,
\end{equation}
where $\Phi$ is an eigenfunction of the Dirichlet Laplacian $\Delta_\Omega$ corresponding to the eigenvalue $\nu_j$.
It is clear that $u_\ell$ satisfies the homogeneous Dirichlet boundary
condition on $\Bbb R\times\partial\Omega$. As $\mu$ meets the condition~\eqref{eq9}, the exponent in~\eqref{test} is an oscillating function of $x\in\Bbb R$. Straightforward calculation shows
that
\begin{equation}\label{s--}
\bigl\|\bigl(\Delta_\Omega-(1+  \lambda)^{-2}(\partial_x+\beta)^2
-\mu\bigr) u_\ell\bigr\|_ {L^2(\Bbb R\times\Omega)}\leq \mathsf c,\quad\|u_\ell\|_{ {H}^2(\Bbb R\times\Omega)}\to\infty
\end{equation}
as $\ell\to +\infty$, where the constant $\mathsf c$ is independent of $\ell$. We extend the functions $u_\ell$ from their supports in $\Pi$ to $\mathcal M$ by zero.

Assume that the estimate~\eqref{coer} is valid. Without loss of generality we can take a rapidly decreasing weight $\mathsf w$, such that the embedding ${\overset{\circ}H}\vphantom{H}^2(\mathcal M) \hookrightarrow L^2(\mathcal M;\mathsf w)$ is compact, and $\|\mathsf w u_\ell\|\leq {\mathsf C}$ uniformly in  $\ell\geq 1$.
Due to stabilization of  ${^\lambda\!\Delta}_\beta$ to $\Delta_\Omega-(1+  \lambda)^{-2}(\partial_x+\beta)^2$ at infinity we have
$$
\|({^\lambda\!\Delta}-\Delta_\Omega+(1+  \lambda)^{-2}(\partial_x+\beta)^2)u_\ell\|\leq\texttt{c}_\ell\|u_\ell\|_{{H}^2(\Bbb R\times\Omega)},
$$
 where $\texttt{c}_\ell\to 0$ as $\ell\to+\infty$. This together with~\eqref{s--}  gives
\begin{equation}\label{s++}
\begin{aligned}
\|({^\lambda\!\Delta}_\beta-  \mu)u_\ell\|\leq C\bigl\|\bigl(-(1+  \lambda)^{-2}(\partial_x+\beta)^2+\Delta_\Omega
-\mu\bigr) u_\ell\bigr\|_ {L^2(\Bbb R\times\Omega)}
\\
+\|({^\lambda\!\Delta}-\Delta_\Omega+(1+  \lambda)^{-2}(\partial_x+\beta)^2)u_\ell\|
\leq C\mathsf c+\texttt{c}_\ell\|u_\ell\|_{{H}^2(\Bbb R\times\Omega)}.
\end{aligned}
\end{equation}
Finally, as a consequence of~\eqref{coer} and~\eqref{s++} we get
\begin{equation}\label{final}
\begin{aligned}
\|u_\ell\|_{{H}^2(\Bbb
R\times\Omega)}
\leq \mathrm C\bigl(\|({^\lambda\!\Delta}-\mu)u_\ell\|+\|\mathsf w u_\ell\|\bigr)
\\
\leq \mathrm C( C\mathsf c
+\texttt{c}_\ell\|u_\ell\|_{ {H}^2(\Bbb
R\times\Omega)}+\mathsf C).
\end{aligned}
\end{equation}
Since $\texttt{c}_\ell\to 0$, the inequalities~\eqref{final} imply that the value $\|u_\ell\|_{ {H}^2(\Bbb
R\times\Omega)}$ remains bounded as $\ell\to+\infty$. This contradicts~\eqref{s--}. Thus the sequence $\{u_\ell\}_{\ell=1}^\infty$ violates the estimate~\eqref{coer}. The necessity is proven.
\qed\end{pf}

\begin{cor}\label{c1}
All eigenvalues of the  Laplacian ${\Delta}$ on a manifold with an asymptotically cylindrical end are of finite multiplicity.
\end{cor}
\begin{pf} By Proposition~\ref{ess} for every $\mu$ there exists $\beta>0$, such that  the  operator ${^0\!\Delta}_\beta-\mu: \mathbf D({^0\!\Delta}_\beta)\to L^2(\mathcal M)$ is Fredholm.   From $\Psi\in \ker ({\Delta}-\mu )$ it follows that $e^{-\beta\mathsf s}\Psi\in \ker({^0\!\Delta}_\beta-\mu)$. As a consequence, $$\dim \ker ({\Delta}-\mu )\leq \dim\ker({^0\!\Delta}_\beta-\mu)<\infty.$$
\qed\end{pf}

\section{Exponential decay of the non-threshold eigenfunctions}\label{s6}

Observe that  the set $\Bbb C\setminus\sigma_{ess}({^\lambda\!\Delta})$ is simply connected. A standard argument based on the analytic Fredholm theory shows that the  spectrum $\sigma({^\lambda\!\Delta})$ of the m-sectorial operator ${^\lambda\!\Delta}$ is the union of the essential spectrum $\sigma_{ess}({^\lambda\!\Delta})$ and the discrete spectrum $\sigma_d({^\lambda\!\Delta})$, e.g.~\cite{KalvinII,Simon Reed iv},~\cite[Appendix]{KozlovMaz`ya}.

As in the theory of $N$-body Schr\"{o}dinger operators, see e.g.~\cite[Chapter XII.11]{Simon Reed iv} and references therein, it is much easier to prove exponential decay of the eigenfunctions corresponding to the isolated eigenvalues. In our case this does not require any assumptions on the analytic regularity of the metric $\mathsf g$, and we can consider a general manifold with an asymptotically cylindrical end in the sense of Definition~\ref{ace}. However, only the Dirichlet Laplacian may have isolated eigenvalues below the first threshold $\nu_1>0$. For the Neumann Laplacian and for the Laplacian on a manifold without boundary the first threshold $\nu_1$ is zero and the eigenvalues are embedded into the essential spectrum $\sigma_{ess}(\Delta)=[0,\infty)$.

In the next lemma we study  the eigenfunctions of $\sigma_d({^\lambda\!\Delta})$. In particular, this lemma implies exponential decay of the eigenfunctions corresponding to the isolated eigenvalues of the Dirichlet Laplacian.

\begin{lem}\label{lemPW1} Assume that $(\mathcal M,\mathsf g)$ is a manifold with an axial analytic asymptotically cylindrical end, and  $\lambda\in\mathcal D_\alpha$ is fixed. Then for  $\mu\in\sigma_{d}({^\lambda\!\Delta})$ and $\Psi\in\ker({^\lambda\!\Delta}-\mu)$ we have $e^{-\beta\mathsf s}\Psi\in \ker({^\lambda\!\Delta}_\beta-\mu)\subset L^2(\mathcal M)$ with some $\beta<0$. In other words, the eigenfunctions corresponding to an isolated eigenvalue of ${^\lambda\!\Delta}$ are exponentially decaying at infinity in the mean.

For a general manifold $(\mathcal M,\mathsf g)$ with an asymptotically cylindrical end the assertion remains valid for $\lambda=0$.
\end{lem}
\begin{pf} Recall from~\cite{Kato,Simon Reed iv} that the family of m-sectorial operators $\Bbb C\ni\beta\mapsto{^\lambda\!\Delta}_\beta$ is said to be analytic of type (B), if for any $\beta\in\Bbb C$ the sectorial form
$\mathsf q^\beta_\lambda$ is densely defined and closed, its domain $\mathbf D(\mathsf q)$ is independent of  $\beta$, and the function $\Bbb C\ni\beta\mapsto\mathsf q^\beta_\lambda[u,u]$ is  analytic for any  $u\in\mathbf D(\mathsf q)$.
Thus the family $\Bbb C\ni\beta\mapsto{^\lambda\!\Delta}_\beta$ is analytic of type (B) by Proposition~\ref{q beta}.
As is known~\cite{Kato,Simon Reed iv}, this implies that for some $\epsilon>0$ there exist functions $\mu_1,\dots,\mu_k$ in the disk $\{\beta\in\Bbb C:|\beta|<\epsilon\}$ with at worst algebraic branching point at $\beta=0$, such that the spectrum of ${^\lambda\!\Delta}_\beta$ in a small neighborhood $\mathcal O$ of $\mu\in\sigma_{d}({^\lambda\!\Delta}_0)$ consists of the isolated eigenvalues $\mu_1(\beta),\dots,\mu_k(\beta)$. Moreover, in the same disk $|\beta|<\epsilon$ there is the analytic function
$$
\beta\mapsto\mathsf P(\beta)=\frac 1 {2\pi i}\oint_{\partial \mathcal O} ({^\lambda\!\Delta}_\beta-\zeta)^{-1}\,d\zeta,
$$
whose values are the projections onto the generalized eigenspace of ${^\lambda\!\Delta}_\beta$ associated with the eigenvalues $\mu_1(\beta),\dots,\mu_k(\beta)$.
For all $\beta\in i\Bbb R$ and $u\in\mathbf D({^\lambda\!\Delta})$ we have $e^{\beta \mathsf s}\,{^\lambda\!\Delta}_\beta\, e^{-\beta\mathsf s}u={^\lambda\!\Delta}u$. Therefore  $\mu_j(\beta)=\cdots=\mu_k(\beta)=\mu$ and
\begin{equation}\label{proj}
\mathsf P(\beta)e^{-\beta\mathsf s}f=e^{-\beta\mathsf s}\mathsf P(0)f,\quad  f\in C_0^\infty(\mathcal M),
\end{equation}
for all $\beta\in i\Bbb R$. By analyticity these equalities  extend to the disk $|\beta|<\epsilon$. The set $C_0^\infty(\mathcal M)$ is dense in $L^2(\mathcal M)$, and the range of $\mathsf P(0)$ is finite dimensional. Hence for any $\Psi\in \ker({^\lambda\!\Delta}-\mu)$ we have $\Psi=\mathsf P(0)f$ with some $f\in C_0^\infty(\mathcal M)$. Due to~\eqref{proj} the function $i\Bbb R\ni\beta \mapsto e^{-\beta\mathsf s}\Psi\in L^2(\mathcal M)$ extends by analyticity to the disk $|\beta|<\epsilon$. Clearly, $e^{-\beta\mathsf s}\Psi\in \ker({^\lambda\!\Delta}_\beta-\mu)$.

It can be shown that the equality~\eqref{proj} and the inclusions $\mu\in\sigma_d({^\lambda\!\Delta}_\beta)$, $e^{-\beta\mathsf s}\Psi\in \ker({^\lambda\!\Delta}_\beta-\mu)$ remain valid as $\beta$ varies in a neighborhood of zero so that  the parabolas of $\sigma_{ess}({^\lambda\!\Delta}_\beta)$ do not cover the point $\mu$, cf. Proposition~\ref{ess} and Fig.~\ref{fig5}. Let us also note that an independent proof of the assertion can be obtained by methods of the asymptotic theory in~\cite{MP2}, see also~\cite[Chapter 5]{KozlovMazyaRossmann}.
\qed\end{pf}

As is well-known e.g.~\cite[Chapter XII.11]{Simon Reed iv}, for $N$-body Schr\"{o}dinger operators with dilation analytic potentials it is also possible to prove exponential decay of all non-threshold eigenfunctions. In the next theorem we show that a similar argument, based on the complex scaling and the Phragm\'{e}n-Lindel\"{o}f principle, allows to prove exponential decay of the non-threshold eigenfunctions of the Laplacian on a manifold with an axial analytic asymptotically cylindrical end.
\begin{thm}\label{exdecay1} Let $\Psi$ be a non-threshold eigenfunction of the Laplacian $\Delta$ on a manifold with an axial analytic asymptotically cylindrical end; i.e. $\Delta\Psi=\mu\Psi\in L^2(\mathcal M)$ with $\mu\in\sigma(\Delta)\setminus\{\nu_j\}_{j=1}^\infty$. Then the  estimate~\eqref{pointwise}
holds for some $\gamma<0$  and an independent of $x$ constant $C$.
\end{thm}

The proof of theorem is preceded by the following lemma.

\begin{lem}\label{lemPW2}  Let the assumptions of Theorem~\ref{exdecay1} be fulfilled. Then the function $\lambda\mapsto \Psi\circ\varkappa_\lambda\in L^2(\mathcal M)$ extends by analyticity from real to all $\lambda$ in the disk $\mathcal D_\alpha$. Moreover, $\Psi\circ\varkappa_\lambda\in\ker({^\lambda\!\Delta}-\mu)$ for all $\lambda\in\mathcal D_\alpha$, and    $\mu\in\sigma_{d}({^\lambda\!\Delta})$ for all non-real $\lambda\in\mathcal D_\alpha$. (Recall that $\varkappa_\lambda$ with $\lambda\in\Bbb R\cap\mathcal D_\alpha$ is the selfdiffeomorphism of $\mathcal M$ defined in Section~\ref{sCS}, and $\alpha<\pi/4$ is an angle for which the conditions of Definition~\ref{ACE} are fulfilled.)
\end{lem}

\begin{pf*}{PROOF of Lemma~\ref{lemPW2}} As  preliminaries to the proof we  briefly recall a construction from~\cite{KalvinII}.

For $\lambda\in\mathcal D_\alpha\cap\Bbb R$ the operator ${^\lambda\!\Delta}$ is the Laplacian on $(\mathcal M,\varkappa_\lambda^*\mathsf g)$, and
the Riemannian geometry gives the identity
\begin{equation}\label{asd}
(\Delta -\zeta)u=\bigl(({^\lambda\!\Delta}-\zeta)(u\circ\varkappa_\lambda)\bigr)\circ\varkappa_\lambda^{-1}\quad \forall u\in \mathbf C({^0\!\Delta}),
\end{equation}
 where $u\circ\varkappa_\lambda$ is in the core $\mathbf C({^\lambda\!\Delta})$ introduced in Definition~\ref{d1}. Let $\zeta<0$ be outside of the sector of the m-sectorial operators ${^\lambda\!\Delta}$, $\lambda\in\mathcal D_\alpha$. Then  the resolvent $({^\lambda\!\Delta}-\zeta)^{-1}$ is an analytic function of $\lambda\in\mathcal D_\alpha$, and we can rewrite~\eqref{asd} in the form
\begin{equation}\label{++}
(\Delta -\zeta)^{-1}F=\bigl(({^\lambda\!\Delta}-\zeta)^{-1}(F\circ\varkappa_\lambda)\bigr)\circ\varkappa_\lambda^{-1},\quad \lambda\in\mathcal D_\alpha\cap\Bbb R,
\end{equation}
where $F$ is in the dense in $L^2(\mathcal M)$ subset $\{F=(\Delta-\zeta)u: u\in \mathbf C({^0\!\Delta})\}$. The equality~\eqref{++} extends by continuity to all $F\in L^2(\mathcal M)$. Taking the inner product of~\eqref{++} with $G\in L^2(\mathcal M)$, we obtain the identity
\begin{equation}\label{h6}
\bigl((\Delta-\zeta)^{-1}F,G
\bigr)=\bigl(({^\lambda\!\Delta}-\zeta)^{-1}(F\circ\varkappa_\lambda),G\circ\varkappa_{\overline{\lambda}}
\bigr)_\lambda
\end{equation}
for all $\lambda\in \mathcal D_\alpha\cap\mathbb R$. The identity~\eqref{h6} cannot be extended by analyticity to all $\lambda\in\mathcal D_\alpha$ for arbitrary  $F,G\in L^2(\mathcal M)$. However it can be done for all $F$ and $G$ in some subset $\mathcal A\subset L^2(\mathcal M)$ of analytic vectors.

In order to introduce the set $\mathcal A$, consider  the algebra $\mathscr E$ of all entire functions $\mathbb C\ni z\mapsto f(z,\cdot)\in C^\infty(\Omega)$, such that
in any sector $|\Im z|\leq (1-\epsilon) \Re z$ with $\epsilon>0$  the
value $\|f(z,\cdot)\|_{L^2(\Omega)}$ decays faster than any  inverse power of $\Re
z$  as $\Re z\to+\infty$. By definition a function $F\in L^2(\mathcal M)$ is in the subset $\mathcal A$ of analytic vectors, if $F(x,\mathrm y)=f(x,\mathrm y)$  for some $f\in\mathscr E$ and all $(x,\mathrm y)\in\Pi$. For $F\in\mathcal A$ and $\lambda\in\mathbb C$  we define the function  $F\circ\varkappa_\lambda$ on $\mathcal M$, such that  $F\circ\varkappa_\lambda\equiv F$ on $\mathcal M\setminus\Pi$, and
\begin{equation}\label{eq}
F\circ\varkappa_\lambda(x,\mathrm y)=f(x+\lambda \mathsf s_R(x),\mathrm y),\quad (x,\mathrm y)\in\Pi.
\end{equation}
Here $f(x+\lambda \mathsf s_R(x),\cdot)$ is the value of the corresponding to $F$ entire function $f\in\mathscr E$ at the point $z=x+\lambda \mathsf s_R(x)$, and $\varkappa_\lambda(x,\mathrm y)=(x+\lambda \mathsf s_R(x),\mathrm y)$ is the complex scaling in $\Pi$. By \cite[Lemma~7.1]{KalvinII} we have:
\begin{itemize}
\item[i.] For any $F\in\mathcal A$, $\mathcal D_\alpha\ni\lambda\mapsto  F\circ\varkappa_\lambda$ is an $L^2(\mathcal M)$-valued analytic function;
\item[ii.] For any $\lambda\in\mathcal D_\alpha$ the image
$\varkappa_\lambda[\mathcal A]=\{F\circ\varkappa_\lambda: F\in\mathcal A\}$
of $\mathcal A$ under $\varkappa_\lambda$ is dense in the space $L^2(\mathcal M)$.
\end{itemize}

For $F,G\in\mathcal A$ the equality~\eqref{h6} extends by analyticity to all $\lambda\in\mathcal D_\alpha$. Let $\mu$ be the same as in the assertion of the lemma. Then $\mu\notin\sigma_{ess}({^\lambda\!\Delta})$ for all non-real $\lambda\in\mathcal D_\alpha$ by Proposition~\ref{ess}.  By applying  the Aguilar-Balslev-Combes argument to the equality~\eqref{h6}, one can see that for any $\lambda\in\mathcal D_\alpha\setminus\Bbb R$ the resolvent $({^\lambda\!\Delta}-\zeta)^{-1}$ is an analytic function of $\zeta$ in a small complex neighborhood of $\mu$, except for the point $\mu\in\sigma_d({^\lambda\!\Delta})\cap\Bbb R$ itself, where the resolvent $({^\lambda\!\Delta}-\zeta)^{-1}$ has a simple pole;  for details we refer to~\cite{KalvinII}. Now the preliminaries are complete, and we are in position to prove the assertion.

Let $\eta\in\mathcal D_\alpha\cap\Bbb R$. Then the Laplacian ${^\eta\!\Delta}$ and the projection
$$
\mathsf P(\eta)=\operatorname{s-lim}\limits_{\epsilon\downarrow 0} i\epsilon ({^\eta\!\Delta}-\mu+i\epsilon)^{-1}
$$
onto its eigenspace  corresponding to the eigenvalue $\mu$ are selfadjoint with respect to the inner product $(\cdot,\cdot)_\eta$ in $L^2(\mathcal M)$. If $\lambda\in\mathcal D_\alpha\setminus\Bbb R$, then we define the projection $\mathsf P(\lambda)$ onto the eigenspace of the non-selfadjoint operator ${^\lambda\!\Delta}$ associated with the eigenvalue $\mu\in\sigma_d({^\lambda\!\Delta})\cap\Bbb R$ as the Riesz projection; i.e. as the first order residue of the resolvent $({^\lambda\!\Delta}-\zeta)^{-1}$ at the simple pole $\mu$. The resolvent $({^\lambda\!\Delta}-\zeta)^{-1}$ is an analytic function of two variables on the set $\{(\lambda,\zeta):\lambda\in\mathcal D_\alpha,\mu\in\Bbb C\setminus\sigma({^\lambda\!\Delta})\}$, e.g.~\cite[Theorem XII.7]{Simon Reed iv}. Therefore the Riesz projection $\mathsf P(\lambda)$ is an analytic function of $\lambda$ on the set $\mathcal  D_\alpha\setminus\Bbb R$.

As a consequence of the equality~\eqref{h6}, for all $F,G\in\mathcal A$ we get
\begin{equation}\label{spare}
\bigl(\mathsf P(0)F,G\bigr)=\bigl(\mathsf P(\eta)(F\circ\varkappa_{\eta}),G\circ\varkappa_{\eta}\bigr)_\eta=\bigl(\mathsf P(\lambda) (F\circ\varkappa_\lambda),G\circ\varkappa_{\overline\lambda}\bigr)_\lambda.
\end{equation}
Recall that  the equality $(\varrho_\lambda \mathcal F,\mathcal G)_\lambda=(\mathcal F,\mathcal G)$ is valid for all $\mathcal F,\mathcal G\in L^2(\mathcal M)$, where $\mathcal D_\alpha\ni\lambda\mapsto\varrho_\lambda\in C^\infty(\mathcal M)$ is an analytic function satisfying~\eqref{0}. Since the sets $\varkappa_\lambda[\mathcal A]$ and $\varkappa_{\overline\lambda}[\mathcal A]$ are dense in $L^2(\mathcal M)$, by the equalities~\eqref{spare} we have
$$
\begin{aligned}
\|& \mathsf P (\lambda)\|=\sup_{F,G\in\mathcal A} \frac{\bigl(\mathsf P(\eta)(F\circ\varkappa_{\eta}),G\circ\varkappa_{\eta}\bigr)_\eta}{\|F\circ\varkappa_\lambda/\sqrt{\varrho_\lambda}\|\|G\circ\varkappa_{\overline{\lambda}}/\sqrt{\varrho_{\overline{\lambda}}}\|}
\\
& \leq \sup_{F,G\in\mathcal A} \frac{\|F\circ\varkappa_\eta/\sqrt{\varrho_\eta}\|\|G\circ\varkappa_\eta/\sqrt{\varrho_\eta}\|}{\|F\circ\varkappa_\lambda/\sqrt{\varrho_\lambda}\|\|G\circ\varkappa_{\overline{\lambda}}/\sqrt{\varrho_{\overline{\lambda}}}\|}\to 1\text{ as } \lambda\to \eta,\lambda\in\mathcal D_\alpha\setminus\Bbb R.
\end{aligned}
$$
Thanks to~\eqref{spare} we also  have
\begin{equation}\label{eq0}
\begin{aligned}
\bigl(\mathsf P&(\lambda)(F\circ\varkappa_\lambda)-\mathsf P(\eta)(F\circ\varkappa_{\eta}),G\circ\varkappa_{\eta}\bigr)_\eta
\\
=&\bigl(\mathsf P(\lambda) (F\circ\varkappa_\lambda),G\circ\varkappa_{\eta}/\varrho_\eta-G\circ\varkappa_{\overline\lambda}/\varrho_{\overline{\lambda}}\bigr)\to 0\text{ as } \lambda\to \eta.
\end{aligned}
\end{equation}
Here the right hand side tends to zero because the norm $\|\mathsf P (\lambda) (F\circ\varkappa_\lambda)\|$ remains bounded, while $G\circ\varkappa_{\overline\lambda}/\varrho_{\overline{\lambda}}$ tends to $G\circ\varkappa_{\eta}/\varrho_\eta$ in  $L^2(\mathcal M)$ as $\lambda\to \eta$. The set $\{G\circ\varkappa_{\eta}:G\in\mathcal A\}$ is dense in $L^2(\mathcal M)$, and hence~\eqref{eq0} implies that $\mathsf P(\lambda) (F\circ\varkappa_\lambda)$ weakly converges to $\mathsf P(\eta)(F\circ\varkappa_\eta)$ as $\lambda\to \eta$, e.g.~\cite[Lemma III.1.31]{Kato}. Therefore  the function $\lambda \mapsto \mathsf P(\lambda) (F\circ\varkappa_\lambda)\in L^2(\mathcal M)$ is weakly (and therefore strongly) analytic in the whole disk $\mathcal D_\alpha$ for any $F\in\mathcal A$. Due to the equality~\eqref{++} we have $\bigl(\mathsf P(0)F\bigr)\circ\varkappa_\lambda=\mathsf P(\lambda)(F\circ\varkappa_\lambda)$ (first for all real, and then by analyticity) for all $\lambda\in\mathcal D_\alpha$.
Since the range $\operatorname{Ran}\mathsf P(0)=\ker (\Delta-\mu)$ is finite dimensional and the set $\mathcal A$ is dense in $L^2(\mathcal M)$, for any $\Psi\in\ker (\Delta-\mu)$ there exists $F\in \mathcal A$, such that $\Psi=\mathsf P(0)F$. The right hand side of the equality $\Psi\circ\varkappa_\lambda=\mathsf P(\lambda)(F\circ\varkappa_\lambda)$ provides the left hand side with an analytic continuation in $\lambda\in\mathcal D_\alpha$. The continuation takes its values in the space $L^2(\mathcal M)$. It remains to note that $\operatorname{Ran} \mathsf P(\lambda)=\ker({^\lambda\!\Delta}-\mu)$ as $\mu\in \sigma_d({^\lambda\!\Delta})\cap\Bbb R$ is a simple pole of the resolvent $({^\lambda\!\Delta}-\zeta)^{-1}$,  e.g.~\cite{Kato}.
\qed\end{pf*}

\begin{pf*}{PROOF of Theorem~\ref{exdecay1}} In the first part of the proof we establish the analyticity of $\Psi$ with respect to $x$ in a complex conical neighborhood of infinity. Then~\eqref{pointwise} follows from a variant of the the Phragm\'{e}n-Lindel\"{o}f principle.

For brevity we denote $\Psi_\lambda=\Psi\circ \varkappa_\lambda$.  By  Lemma~\ref{lemPW2} the function $\mathcal D_\alpha\ni\lambda \mapsto \Psi_\lambda\in L^2(\mathcal M)$  is analytic, and $\Psi_\lambda\in\ker({^\lambda\!\Delta}-\mu)$. (Note that $\ker({^\lambda\!\Delta}-\mu)\subset C^\infty(\mathcal M)$ by usual results on local regularity of solutions to elliptic problems.) We have
$$
\mathsf q_\lambda  [\Psi_\lambda, v]-\mu(\Psi_\lambda,v)_\lambda=0\quad \forall v\in\mathbf D(\mathsf q).
$$
This together with Proposition~\ref{p1}.ii and~\eqref{3} gives
$$
\|\Psi_\lambda \|^2_{\mathbf D(\mathsf q)}=(d \Psi_\lambda, d\Psi_\lambda)+\|\Psi_\lambda\|^2\leq ( b^{-1}(c\mu+a)+1)\|\Psi_\lambda\|^2.
$$
 As is known~\cite[Chapter VII.4]{Kato}, results of Proposition~\ref{p1} also imply that
$$
|\mathsf q_\lambda[u,v]-\mathsf q_\varsigma[u,v]|\leq C_{\lambda,\varsigma} \|u\|_{\mathbf D(\mathsf q)}\|v\|_{\mathbf D(\mathsf q)},
$$
where the constant $C_{\lambda,\varsigma}$ is bounded uniformly in $u,v\in\mathbf D(\mathsf q)$ and $\lambda,\varsigma\in\mathcal D_\alpha$; moreover,  $C_{\lambda,\varsigma}\to 0$  as $|\lambda-\varsigma|\to 0$.
As a consequence we obtain
$$
\begin{aligned}
b\|\Psi_\lambda-\Psi_\varsigma\|^2_{\mathbf D(\mathsf q)}- (a+b+c\mu)\|\Psi_\lambda-\Psi_\varsigma\|^2
\\\leq \Re\mathsf q_\lambda[\Psi_\lambda-\Psi_\varsigma,\Psi_\lambda-\Psi_\varsigma]-\mu\Re(\Psi_\lambda-\Psi_\varsigma,\Psi_\lambda-\Psi_\varsigma)_\lambda
\\
\leq \bigl|\mathsf q_\lambda[\Psi_\varsigma,\Psi_\lambda -\Psi_\varsigma]-\mu(\Psi_\varsigma,\Psi_\lambda-\Psi_\varsigma)_\lambda-\mathsf q_\varsigma[\Psi_\varsigma,\Psi_\lambda-\Psi_\varsigma]+\mu(\Psi_\varsigma,\Psi_\lambda-\Psi_\varsigma)_\varsigma\bigr|
\\
\leq \bigl|\mathsf q_\lambda[\Psi_\varsigma,\Psi_\lambda -\Psi_\varsigma]-\mathsf q_\varsigma[\Psi_\varsigma,\Psi_\lambda-\Psi_\varsigma]\bigr|
+\mu\bigl|\bigl((\varrho_\lambda-\varrho_\varsigma)\Psi_\varsigma,\Psi_\lambda-\Psi_\varsigma\bigr)\bigr|\\
\leq \mathrm C_{\lambda,\varsigma} \|\Psi_\varsigma\|
\|\Psi_\lambda-\Psi_\varsigma\|_{\mathbf D(\mathsf q)},
\end{aligned}
$$
where $\mathrm C_{\lambda,\varsigma}$ is uniformly bounded, and  $\mathrm C_{\lambda,\varsigma}\to 0$ as $|\lambda-\varsigma|\to 0$. By these estimates the analyticity  of $\mathcal D_\alpha\ni\lambda\mapsto \Psi_\lambda\in L^2(\mathcal M)$ leads to the continuity of the function $\mathcal D_\alpha\ni\lambda\mapsto \Psi_\lambda\in\mathbf D(\mathsf q)$. Then by the Morera's theorem the  function $\mathcal D_\alpha\ni\lambda\mapsto \Psi_\lambda\in\mathbf D(\mathsf q)$ is analytic. The graph norm of $\Delta^{1/2}$ is an equivalent norm in $\mathbf D(\mathsf q)$ e.g.~\cite{Kato}, and therefore usual results on traces of functions in the Sobolev space $H^1_{loc}(\mathcal M)\supseteq\mathbf D(\mathsf q)$ apply. In particular, for any fixed $x\in \Bbb R_+$ we have
$$
\|u(x)\|_{L^2(\Omega)}\leq\|u(x)\|_{H^{1/2}(\Omega)}\leq c\|u\|_{\mathbf D(\mathsf q)},
$$
 where  $c$ is independent of $u\in \mathbf D(\mathsf q)$. Hence for any $x\in \Bbb R_+$ the function  $\mathcal D_\alpha\ni\lambda\mapsto \Psi_\lambda(x)\in L^2(\Omega)$  is  analytic.

Consider $\Psi_\lambda(x)=\Psi(x+\lambda\mathsf s_R(x))$  as an analytic $L^2(\Omega)$-valued function of  $z=x+\lambda\mathsf s_R(x)$ in the  complex conical neighborhood
$$
\mathcal S_\epsilon=\{x+\lambda\mathsf s_R(x)\in\Bbb C: x\in\Bbb R_+, \mathsf s_R(x)\neq 0, |\lambda|\leq\epsilon\}, \quad\epsilon<\sin\alpha,
$$
of infinity.
 The function $\Psi$ is uniformly bounded on $\mathcal S_\epsilon$  in  the sense that
\begin{equation}\label{aux1}
\int_{\mathfrak L_\lambda^R}\|\Psi(z)\|_{L^2(\Omega)}^2\,|dz|=\int_0^\infty \|\Psi_\lambda(x)\|^2_{L^2(\Omega)} |1+\lambda\mathsf s'_R(x)|\, dx\leq 2\|\Psi_\lambda\|^2\leq C,
\end{equation}
where $\mathfrak L^R_\lambda$ is the  curve depicted on Fig.~\ref{fig+}, and $|\lambda|\leq \epsilon$.  By Lemmas~\ref{lemPW2} and~\ref{lemPW1} for any non-real $\lambda\in\mathcal D_\alpha$ we have  $\mu\in\sigma_{d}({^\lambda\!\Delta})$ and $e^{-\beta(\lambda)\mathsf s}\Psi_\lambda\in L^2(\mathcal M)$  with  some negative $\beta(\lambda)$. Therefore the function $\mathcal S_\epsilon\ni z\mapsto \Psi(z)\in L^2(\Omega)$ is exponentially decaying outside of the half-axis $\Bbb R_+$ in the sense that
\begin{equation}\label{aux2}
\int_{\mathfrak L_\lambda^R}\|e^{-\beta(\lambda) z/2}\Psi(z)\|_{L^2(\Omega)}^2\,|dz|\leq C(\lambda)<\infty,\quad\beta(\lambda)<0, |\lambda|\leq\epsilon,\lambda\notin\Bbb R.
\end{equation}
 In the remaining part of the proof we derive a variant of the Phragm\'{e}n-Lindel\"{o}f principle, which says that  any analytic function $\mathcal S_\epsilon\ni z\mapsto \Psi(z)\in L^2(\Omega)$ satisfying~\eqref{aux1} and~\eqref{aux2} also meets the estimate~\eqref{pointwise} with some $\gamma<0$.

The uniform estimate~\eqref{aux1} necessitates existence of a sequence of positive numbers $\{T_\ell\}_{\ell=1}^\infty$, $T_\ell\to\infty$ as $\ell\to\infty$, such that for a given $\gamma<0$ and any $\delta>0$ the integral
$$
\int_{\{z\in \mathcal  S_\epsilon: z=T_\ell+\lambda \mathsf s_R(T_\ell), |\lambda|\leq \epsilon, |1+\lambda|=1\}} \frac {e^{-\delta z^2-\gamma z}\Psi(z)}{x-z}\, dz, \quad x\in \Bbb R_+,
$$
tends to zero in $L^2(\Omega)$ as $\ell\to\infty$. As a consequence, the contour of the Cauchy integral in the equality
$$
e^{-\delta x^2-\gamma x}\Psi(x)=\frac 1 {2\pi i}\oint_{|x-z|=\wp(x)} \frac {e^{-\delta z^2-\gamma z}\Psi(z)}{x-z}\, dz,\quad x\in\mathcal S_\epsilon\cap\Bbb R_+,
$$
(where $\wp(x)$ is so small that the contour of integration lies in $\mathcal S_\epsilon$)
can be deformed so that
\begin{equation}\label{aux3}
e^{-\delta x^2-\gamma x}\Psi(x)=\frac 1 {2\pi i}\Bigl(\int_{\mathfrak L^R_{i\epsilon}}-\int_{\mathfrak L^R_{-i\epsilon}}\Bigr) \frac {e^{-\delta z^2-\gamma z}\Psi(z)}{x-z}\, dz.
\end{equation}
Here both integrals are absolutely convergent in $L^2(\Omega)$ because of~\eqref{aux1}. Moreover, thanks to~\eqref{aux2}, for $0>\gamma>\max\{\beta(i\epsilon),\beta(-i\epsilon)\}/2$ these integrals   are bounded in $L^2(\Omega)$ uniformly in $\delta>0$ and large $x>0$. Hence the same is true for the left hand side of the equality~\eqref{aux3}. This immediately leads to~\eqref{pointwise}.
\qed\end{pf*}

\section{Refined exponential decay of the non-threshold eigenfunctions and accumulation of eigenvalues}\label{s7}
In this  section we consider the Laplacian on a general manifold $(\mathcal M,\mathsf g)$ with an asymptotically cylindrical end $(\Pi,\mathsf g\!\upharpoonright_{\Pi})$ in the sense of Definition~\ref{ace}. The axial analyticity of the end is not assumed. We study the non-threshold exponentially decaying eigenfunctions and accumulation of the corresponding eigenvalues. Let us stress that in the case of an axial analytic asymptotically cylindrical end  $(\Pi,\mathsf g\!\upharpoonright_{\Pi})$  all non-threshold eigenfunctions of the Laplacian are of some exponential decay by Theorem~\ref{exdecay1}.

\begin{thm}[Refined exponential decay]\label{decay} Let $\Psi$ be a non-threshold exponentially decaying eigenfunction of the Laplacian on a manifold $(\mathcal M,\mathsf g)$ with an asymptotically cylindrical end; i.e. $\Delta\Psi=\mu\Psi$ with $\mu\in\sigma(\Delta)\setminus\{\nu_j\}_{j=1}^\infty$, and the estimate~\eqref{pointwise} holds for some~$\gamma<0$.
Then the estimate~\eqref{pointwise} holds for any negative $\gamma >-\min_{j:\nu_j>\mu}\sqrt{\nu_j-\mu}$.
\end{thm}
\begin{pf} In the proof we use methods of the asymptotic theory~\cite{KozlovMaz`ya,KozlovMazyaRossmann,MP2}.

Let us consider the case of the Dirichlet Laplacian on a manifold $(\mathcal M,\mathsf g)$ with non-compact boundary. The case $\partial\mathcal M=\varnothing$ is similar, while in the case of the Neumann Laplacian some changes in our argument are needed, which are outlined in the end of the proof.

Take  a  cutoff function  $\chi\in C^\infty(\Bbb R)$, such that $0\leq \chi(x)\leq 1$,  $\chi(x)=0$ for $x\leq 1$, and $\chi(x)=1$ for $x\geq 2$. Consider the infinite cylinder $\Bbb R\times\Omega$ endowed with the compound metric  $(1-\chi_T)\bigl( dx\otimes dx+\mathfrak h\bigr)+\chi_T \mathsf g$, where $\chi_T(x)=\chi(x-T)$ and $T$ is a sufficiently large positive number. The compound metric is well-defined,  because the metric $\mathsf g$ stabilizes to the product metric $dx\otimes dx+\mathfrak h$ at infinity, see Definition~\ref{ace}. Let $\Delta^T$ be the Laplacian induced  on $\Bbb R\times\Omega$ by the compound metric. Introduce the conjugated Laplacian $\Delta^T_\beta=e^{-\beta x}\Delta^T e^{\beta x}$, where $e^{\pm\beta x}$ stands for the operator of multiplication  by the exponent. Let ${\overset{\circ}H}\vphantom{H}^2(\Bbb R\times\Omega)$ and $L^2(\Bbb R\times\Omega)$ be the spaces introduced in the proof of Proposition~\ref{ess}. The Dirichlet Laplacian  $\Delta^T_\beta:{\overset{\circ}H}\vphantom{H}^2(\Bbb R\times\Omega)\to L^2(\Bbb R\times\Omega)$ possesses the properties:
\begin{itemize}
\item[1.] $\Delta^T_\beta$ coincides with ${^0\!\Delta_\beta}$ on $(T+1,\infty)\times\Omega\subset \Pi$.
\item[2.] $\Delta^T_\beta$ coincides with $\Delta_\Omega-(\partial_x+\beta)^2$ on $(-\infty,T)\times\Omega$.
\item[3.] As a consequence of stabilization of $\mathsf g$ to $dx\otimes dx+\mathfrak h$ at infinity the estimate
\begin{equation}\label{con}
\|(\Delta_\beta^T-\Delta_\Omega+(\partial_x+\beta)^2)u\|_{L^2(\Bbb R\times\Omega)}\leq C(T)\|u\|_{H^2(\Bbb R\times\Omega)}
\end{equation}
is valid, where the constant $C(T)$ is independent of $u\in {\overset{\circ}H}\vphantom{H}^2(\Bbb R\times\Omega)$ and $\beta\in\mathcal K$, where $\mathcal K$ is a compact subset of $\Bbb C$. Moreover,  $C(T)\to 0$ as $T\to+\infty$, cf.~\eqref{LR} and~\eqref{stab}.
\end{itemize}
Let us show that for arbitrarily small $\delta>0$ there exists a sufficiently large $T>0$, such that the resolvent
\begin{equation}\label{res}
(\Delta^T_\beta-\mu)^{-1}\in \mathscr B(L^2(\Bbb R\times\Omega), {\overset{\circ}H}\vphantom{H}^2(\Bbb R\times\Omega))
\end{equation}
 is an analytic function of $\beta$ on the compact subset
$$
\mathcal K_\delta=\{\beta\in \Bbb C: -\delta\geq\Re\beta\geq\delta-\min_{j:\nu_j>\mu}\sqrt{\nu_j-\mu},\ |\Im\beta|\leq\delta\}
$$
of the complex plane; i.e. it is analytic in a small open neighborhood of $\mathcal K_\delta$.
Recall from the proof of Proposition~\ref{ess}  that the operator~\eqref{opo} yields an isomorphism, if $\beta\in\Bbb C$ does not satisfy  the equality~\eqref{eq9} with any $j\in\Bbb N$ and any $\xi\in\Bbb R$. Thus the operator~\eqref{opo} with $\lambda=0$ and $\mu\in\Bbb R\setminus\{\nu_j\}_{j=1}^\infty$ is an analytic Fredholm operator function of $\beta\in\mathcal K_\delta$, which is invertible for all $\beta\in\mathcal K_\delta$.   The analytic Fredholm theory immediately implies that the inverse operator
$$
(\Delta_\Omega-(\partial_x+\beta)^2-\mu)^{-1}:L^2(\Bbb R\times\Omega)\to {\overset{\circ}H}\vphantom{H}^2(\Bbb R\times\Omega)
$$
is uniformly bounded (analytic) function of $\beta\in\mathcal K_\delta$. From this together with~\eqref{con} we conclude that  for all sufficiently large $T>T(\delta)$ and all $\beta\in\mathcal K_\delta$ the operator norm of the composition
$$
\Lambda=(\Delta_\Omega-(\partial_x+\beta)^2-\Delta_\beta^T) (\Delta_\Omega-(\partial_x+\beta)^2-\mu)^{-1}\in\mathscr B(L^2(\Bbb R\times\Omega))
$$
 is less than one. Therefore for all $\beta\in\mathcal K_\delta$ we have
$$
(\Delta_\beta^T-\mu)^{-1}= (\Delta_\Omega-(\partial_x+\beta)^2-\mu)^{-1}\sum_{j=0}^\infty \Lambda^j,
$$
where the series $\sum_{j=0}^\infty \Lambda^j$ converges  in the space $\mathscr B(L^2(\Bbb R\times\Omega))$.
As $\Delta_\beta^T$ depends on $\beta\in\Bbb C$ analytically, and the inclusion~\eqref{res} is valid for all $\beta\in\mathcal K_\delta$, the resolvent~\eqref{res} is an analytic function of $\beta\in\mathcal K_\delta$.

Let $\Psi$ meet the assumptions of theorem. Thanks to~\eqref{pointwise} for any $\beta\in\Bbb C$ with $\Re\beta>\gamma$ we have $e^{-\beta\mathsf s}\Psi\in \ker(\Delta_\beta-\mu)$.  Then for $\varphi_T=\chi_{T+2}$ we obtain
$$(\Delta_\beta^T-\mu)\varphi_T e^{-\beta\mathsf s}\Psi
=e^{-\beta\mathsf s}[ \Delta,\varphi_T] \Psi.$$
(Here the right hand side of the equality reads as follows: the cutoff function $\varphi_T$ is extended from its support in $\Pi$ to the  manifold  $\mathcal M$ by zero, then  the commutator is well-defined, and we extend the right hand side of the equality from its support in $\Pi$ to $\Bbb R\times\Omega$ by zero.) Consequently,  for all $\beta\in\mathcal K_\delta$, $\Re\beta>\gamma$, the equality
\begin{equation}\label{star}
\varphi_T e^{-\beta\mathsf s}\Psi=
(\Delta_\beta^T-\mu)^{-1}e^{-\beta\mathsf s}[ \Delta,\varphi_T]\Psi
\end{equation}
is valid. The function $[ \Delta,\varphi_T] \Psi\in C^\infty(\Bbb R\times\Omega)$ is compactly supported, and hence $e^{-\beta\mathsf s}[ \Delta,\varphi_T] \Psi$ is an analytic function of $\beta
\in\Bbb C$ with values in $L^2(\Bbb R\times\Omega)$. Since the resolvent~\eqref{res} is analytic in $\beta\in\mathcal K_\delta$, the right hand side of~\eqref{star} is an analytic function of $\beta\in\mathcal K_\delta$ with values in  the space ${\overset{\circ}H}\vphantom{H}^2(\Bbb R\times\Omega)$. Hence the function $\mathcal K_\delta\ni\beta\mapsto \varphi_T e^{-\beta\mathsf s}\Psi\in {\overset{\circ}H}\vphantom{H}^2(\Bbb R\times\Omega)$ is analytic.  As $\delta$ is  arbitrarily small, this establishes the inclusion $\varphi_T e^{-\beta\mathsf s}\Psi\in {\overset{\circ}H}\vphantom{H}^2(\Bbb R\times\Omega)$, where $\Re\beta>-\min_{j:\nu_j>\mu}\sqrt{\nu_j-\mu}$.
Now we can conclude that $\Psi$ meets the pointwise estimate~\eqref{pointwise} with an independent of $x$ constant $C$ and $\gamma>-\min_{j:\nu_j>\mu}\sqrt{\nu_j-\mu}$.

Indeed, let $\hat u(\tau)=\mathcal F_{x\mapsto \tau} u(x)$ be the Fourier transform of  $u=\varphi_T e^{-\beta\mathsf s}\Psi$. Then  $\int_{\Bbb R} (1+\tau)^4 \|\hat u(\tau)\|_{L^2(\Omega)}^2\,d\tau \leq c \|u\|_{H^2(\Bbb R\times\Omega)}$ e.g.~\cite{Lions Magenes}. For
$u(x)=\int_{\Bbb R} e^{i x \tau }\hat u(\tau)\,d\tau$ we deduce the estimates
$$
\begin{aligned}
&\|u(x)\|_{L^2(\Omega)}\leq \int_{\Bbb R}\| e^{i x \tau }\hat u(\tau)\|_{L^2(\Omega)}\,d\tau \leq\left(\int_{\Bbb R} (1+\tau)^{-2}\,d\tau\right)^{1/2}\\
& \times\left( \int_{\Bbb R} (1+\tau)^2 \|\hat u(\tau)\|_{L^2(\Omega)}^2\,d\tau\right)^{1/2}\leq C \|u\|_{H^2(\Bbb R\times\Omega)},
\end{aligned}
$$
where we used the Cauchy-Schwarz inequality. This completes the proof for the case of the Dirichlet Laplacian.

In order to study the case of the Neumann Laplacian,  we consider the continuous operator
\begin{equation}\label{opo1}
\begin{aligned}
\{\Delta_\Omega- &(1+\lambda)^{-2}(\partial_x+\beta)^2 - \mu,\partial_\eta\}:\\
& H^2(\Bbb R\times\Omega)\to L^2(\Bbb R\times\Omega)\times H^{1/2}(\Bbb R\times\partial\Omega)
\end{aligned}
\end{equation}
of the non-homogeneous Neumann problem on $(\Bbb R\times\Omega,dx\otimes dx+\mathfrak h)$. Here  $\partial_\eta$ is the operator of the  Neumann boundary conditions, the space $H^\ell(\Bbb R\times\Omega)$ is introduced as the completion of the set $C_c^\infty(\Bbb R\times\Omega)$ in the norm~\eqref{norm}, and  $H^{1/2}(\Bbb R\times\partial\Omega)$ is the space of traces of the functions in $H^1(\Bbb R\times\Omega)$. Applying the Fourier transform $\mathcal F_{x\mapsto \xi}$ we pass from the operator~\eqref{opo1} to the operator $\{\Delta_\Omega+(1+\lambda)^{-2}(\beta+i\xi)^2-\mu, \partial_\eta\}$ of the non-homogeneous Neumann problem on $(\Omega,\mathfrak h)$. Suppose that  $\mu$ does not satisfy the equality~\eqref{ess} for any $j\in\Bbb N$ and $\xi\in\Bbb R$  or, equivalently, suppose that  for any $\xi\in\Bbb R$ the number $\mu-(1+\lambda)^{-2}(\beta+i\xi)^2$  is not an eigenvalue $\nu_j$ of the Neumann Laplacian $\Delta_\Omega$ on $(\Omega,\mathfrak h)$.  Then a known argument, see e.g.~\cite[Theorem 4.1]{MP2} or \cite[Theorem 5.2.2]{KozlovMazyaRossmann} or \cite[Theorem~2.4.1]{KozlovMaz`ya}, shows that the operator~\eqref{opo1} realizes an isomorphism.

In the same way as before we introduce the compound metric on $\Bbb R\times\Omega$ and consider the corresponding operator $\partial_\nu^T$ of the Neumann boundary condition on $\Bbb R\times\partial\Omega$. The continuous operator
\begin{equation}\label{aux00}
 \begin{aligned}
\{\Delta_\beta^T-\mu,e^{-\beta x}\,\partial^T_\nu\, e^{\beta x}\}:H^2(\Bbb R\times\Omega)\to L^2(\Bbb R\times\Omega)\times H^{1/2}(\Bbb R\times\partial\Omega)
\end{aligned}
\end{equation}
possesses the properties:
\begin{itemize}
\item[1.] It coincides with $\{{^0\!\Delta_\beta}-\mu, e^{-\beta\mathsf s}\partial_\nu e^{\beta\mathsf s}\}$ on $(T+1,\infty)\times\Omega\subset \Pi$.
\item[2.] It coincides with $\{\Delta_\Omega-(\partial_x+\beta)^2-\mu,\partial_\eta\}$ on $(-\infty,T)\times\Omega$.
\item[3.] As a consequence of stabilization of $\mathsf g$ to $dx\otimes dx+\mathfrak h$ at infinity the estimate
\begin{equation*}\label{con1}
\begin{aligned}
\|(\Delta_\beta^T-\Delta_\Omega+(\partial_x+ & \beta)^2)u\|_{L^2(\Bbb R\times\Omega)}
\\
& +\|\partial_\eta u-e^{-\beta x}\,\partial^T_\nu\, e^{\beta x}u\|_{H^{1/2}(\Bbb R\times\partial\Omega)}
\leq C(T)\|u\|_{H^2(\Bbb R\times\Omega)}
\end{aligned}
\end{equation*}
is valid, where the constant $C(T)$ is independent of $u\in {H}^2(\Bbb R\times\Omega)$ and $\beta\in\mathcal K_\delta$; moreover,  $C(T)\to 0$ as $T\to+\infty$.
\end{itemize}
Similarly to the case of the Dirichlet Laplacian one can show that for any $\delta>0$ there exists $T>T(\delta)$, such that the inverse of the operator~\eqref{aux00} is an analytic function of $\beta\in\mathcal K_\delta$. Therefore the equality
$$
\varphi_T e^{-\beta\mathsf s}\Psi=\{\Delta_\beta^T-\mu,e^{-\beta x}\,\partial^T_\nu\, e^{\beta x}\}^{-1} e^{-\beta\mathsf s}\{[ \Delta,\varphi_T]\Psi,  [ \partial_\nu ,\varphi_T]\Psi\}
$$
extends by analyticity from the set $\{\beta\in\mathcal K_\delta: \Re\beta>\gamma\}$ to all $\beta\in\mathcal K_\delta$. This equality is a substitution for~\eqref{star}, it implies the inclusion $\chi e^{-\beta\mathsf s}\Psi\in {H}^2(\Bbb R\times\Omega)$, where $\beta>-\min_{j:\nu_j>\mu}\sqrt{\nu_j-\mu}$.
\qed
\end{pf}
\begin{thm}[Accumulation of eigenvalues]\label{acc}  Let $\{\mu_k\}_{k=1}^\infty$ be a sequence of eigenvalues of the Laplacian $\Delta$ on a manifold $(\mathcal M,\mathsf g)$ with an asymptotically cylindrical end, such that $\mu_k\to\mu<\infty$ as $k\to \infty$, and $\mu_k\neq\mu$. Assume that
to every $\mu_k$ there corresponds an eigenfunction $\Psi_k$ satisfying the estimate~\eqref{pointwise} with some $\gamma=\gamma_k<0$.
Then $\gamma_k\to 0$, $\mu$ is a threshold $\nu_j$ of $\Delta$, and the sequence $\{\mu_k\}_{k=1}^\infty$  accumulates to $\nu_j$ only from below.
\end{thm}

\begin{pf}  Here we combine the ideas used in the proof of Theorem~\ref{decay} with the compactness argument due to Perry~\cite{Perry}.

Let us first show that $\gamma_k\to 0$ as $k\to\infty$. Assume the contrary. Then there exists $\beta<0$, such that $\gamma_k<\beta$ for  all $k>0$. By taking a larger  negative $\beta$ we can always achieve  $\mu\notin\sigma_{ess}({^0\!\Delta_\beta})$, cf. Proposition~\ref{ess} and Fig.~\ref{fig5}.

In the case of the Dirichlet Laplacian, and also in the case  $\partial\mathcal M=\varnothing$, similarly to~\eqref{star} we deduce
$$
\varphi_T e^{-\beta\mathsf s}\Psi_k=(\Delta_\beta^T-\mu)^{-1} (\Delta_\beta-\mu_k+\mu_k-\mu)\varphi_T e^{-\beta\mathsf s}\Psi_k
$$
$$=(\Delta_\beta^T-\mu)^{-1}\bigl(e^{-\beta\mathsf s}[ \Delta,\varphi_T]+(\mu_k-\mu)\varphi_T e^{-\beta\mathsf s}\bigr)\Psi_k.
$$
For all sufficiently large $T$ the resolvent~\eqref{res} is bounded by the argument in the proof of Theorem~\ref{decay}. Since  the metric $\mathsf g$ stabilizes to the product metric at infinity, the norm $\|f\|$ of a function supported in $\Pi$ is equivalent to the norm $\|f\|_{L^2(\Bbb R\times\Omega)}$, and the norm in ${H}^2(\mathcal M)$ is equivalent to the graph norm of the Laplacian, see the proof of Proposition~\ref{ess}. Hence
\begin{equation}\label{V}
\|\varphi_T e^{-\beta\mathsf s}\Psi_k\|^2\leq C(\mu_k^2+1)\|\Psi_k\|^2+C(\mu_k-\mu)^2\|\varphi_T e^{-\beta\mathsf s}\Psi_k\|^2,
\end{equation}
where we extended the functions from their supports to $\mathcal M$ by zero and used the estimates
$$
\|e^{-\beta\mathsf s}[ \Delta,\varphi_T]\Psi_k\|^2\leq c_1\|\Psi_k\|^2_{H^2(\mathcal M)}\leq c_2(\|\Delta \Psi_k\|^2+\| \Psi_k\|^2)=c_2(\mu_k^2+1)\|\Psi_k\|^2.
$$

In the case of the Neumann Laplacian we use the equality
$$
\varphi_T e^{-\beta\mathsf s}\Psi_k=\{\Delta_\beta^T-\mu,e^{-\beta x}\,\partial^T_\nu\, e^{\beta x}\}^{-1} e^{-\beta\mathsf s}\{[ \Delta,\varphi_T]\Psi_k+(\mu_k-\mu)\Psi_k,  [ \partial_\nu ,\varphi_T]\Psi_k\}.
$$
This equality together with boundedness of the inverse of~\eqref{aux00} and the estimates
$$
\begin{aligned}
\|e^{-\beta\mathsf s}[ \Delta,\varphi_T]\Psi_k\|^2+ & \|e^{-\beta\mathsf s}[ \partial_\nu ,\varphi_T]\Psi_k\|^2_{H^{1/2}(\Bbb R\times\partial\Omega)}\leq c_1\|\Psi_k\|^2_{H^2(\mathcal M)}
 \\
& \leq c_2(\|\Delta \Psi_k\|^2+\| \Psi_k\|^2)=c_2(\mu_k^2+1)\|\Psi_k\|^2
\end{aligned}
$$
justifies the uniform in $k$ inequality~\eqref{V} for the Neumann Laplacian.

Clearly, $1-C(\mu_k-\mu)^2>0$ for all sufficiently large $k$, and from~\eqref{V} we get
\begin{equation}\label{123}
\|\varphi_T e^{-\beta\mathsf s}\Psi_k\|^2\leq C(\mu_k^2+1)\bigl(1-C(\mu_k-\mu)^2\bigr)^{-1}\|\Psi_k\|^2\leq Const \|\Psi_k\|^2.
\end{equation}
For  an independent of $\Psi\in L^2(\mathcal M)$ constant $C$ we have $\|(1-\varphi_T)e^{-\beta\mathsf s}\Psi\|\leq C\|\Psi\|$. Thus~\eqref{123} leads to the uniform in $k$ estimate
$\|e^{-\beta\mathsf s}\Psi_k\| \leq \mathrm C\|\Psi_k\|$. This estimate together with the Cauchy-Schwarz inequality gives
$$
\|\Psi_k\|^2\leq \|e^{-\beta\mathsf s}\Psi_k\|\|e^{\beta\mathsf s}\Psi_k\|\leq  \mathrm C\|\Psi_k\|\,\|e^{\beta\mathsf s}\Psi_k\|.
$$
Finally, for all sufficiently large $k$ we obtain the uniform estimate
\begin{equation}\label{*}
\|\Psi_k\|^2_{H^2(\mathcal M)}\leq C(\|\Delta \Psi_k\|^2+\| \Psi_k\|^2)=C(\mu_k^2+1)\|\Psi_k\|^2\leq Const \|e^{\beta\mathsf s}\Psi_k\|^2.
\end{equation}
The Sobolev space $H^2(\mathcal M)$ is compactly embedded into the weighted space $L^2(\mathcal M,e^{\beta\mathsf s})$. Since compactness of the unit ball implies that the space is finite-dimensional, there can only be a finite number of linearly independent normalized eigenfunctions $\Psi_k$ satisfying~\eqref{*}. Since $\mu_k\to \mu\neq \mu_k $ as $k\to\infty$,  we come to a contradiction.

We proved that $\gamma_k\to 0$ as $k\to\infty$. It remains to note that for any subsequence of $\{\mu_k\}_{k=1}^\infty$ that does not accumulate to a threshold $\nu_j$ from below we can use Theorem~\ref{decay} in order to refine the corresponding subsequence of $\{\gamma_k\}_{k=1}^\infty$ so that the elements $\gamma_k$ of the subsequence does not tend to zero as $k\to\infty$. This completes the proof.

Let us remark here that an independent proof of this theorem can be obtained by methods announced in~\cite{KNP} and then developed in~\cite{Kalvine,Plam} and~\cite[Appendix]{PhD}.
\qed\end{pf}

\begin{pf*}{Proof of Theorem~\ref{main}.} The assertions are readily apparent from Theorems~\ref{exdecay1},~\ref{decay},~\ref{acc}, and Corollary~\ref{c1}.
\qed\end{pf*}

\noindent{\it Acknowledgement.} It is a pleasure to thank Werner M\"{u}ller for useful discussions and for his hospitality at the University of Bonn.

\end{document}